\documentclass[11 pt, psamsfonts]{amsart}
\usepackage{amssymb}

\def\ignore#1{\relax}

\def\g{\mathfrak g}
\def\b{\mathfrak b}
\def\h{\mathfrak h}
\def\n{\mathfrak n}

\def\sl{\mathfrak{sl}}
\def\so{\mathfrak{so}}

\def\Vl{{V_\la}}
\def\R{{\mathbb R}}
\def\Z{{\mathbb Z}}

\def\Q{{\mathbb Q}}
\def\C{{\mathbb C}}

\def\la{\lambda}
\def\ep{\epsilon}

\def\r{{\bf r}}

\def\B{{\mathcal B}}

\def\ni{\noindent}

\def\ignore#1{\relax}
\def\om{\omega}

\def\1{{\bf 1}}

\def\End{{\rm End}}

\calclayout

\renewenvironment{proof}[1][\proofname]{\par
  \normalfont
  \topsep6\p@\@plus6\p@ \trivlist
  \item[\hskip\labelsep\bfseries
    #1\@addpunct{.}]\ignorespaces
}{%
  \qed\endtrivlist
}

\def\th@plain{%
  \let\thmhead\thmhead@plain \let\swappedhead\swappedhead@plain
  \thm@preskip.5\baselineskip\@plus.2\baselineskip
                                    \@minus.2\baselineskip
  \thm@postskip\thm@preskip
  \itshape
\renewcommand{\labelenumi}{{(\alph{enumi})\quad}}
                        \renewcommand{\labelenumii}{{(\roman{enumii})\ }}
}
\def\th@definition{%
  \let\thmhead\thmhead@plain \let\swappedhead\swappedhead@plain
  \thm@preskip.5\baselineskip\@plus.2\baselineskip
                                    \@minus.2\baselineskip
  \thm@postskip\thm@preskip
  \upshape
}
\def\th@remark{%
  \thm@headfont{\itshape}
  \let\thmhead\thmhead@plain \let\swappedhead\swappedhead@plain
  \thm@preskip.5\baselineskip\@plus.2\baselineskip
                                    \@minus.2\baselineskip
  \thm@postskip\thm@preskip
  \upshape
}

{\theoremstyle{plain}
\newtheorem{theorem}{Theorem}[section]
}

{\theoremstyle{plain}
\newtheorem{proposition}[theorem]{Proposition}
}

{\theoremstyle{plain}
\newtheorem{corollary}[theorem]{Corollary}
}

{\theoremstyle{plain}
\newtheorem{lemma}[theorem]{Lemma}
}

{\theoremstyle{plain}

}

{\theoremstyle{definition}
\newtheorem{definition}[theorem]{Definition}
}

{\theoremstyle{definition}

}

{\theoremstyle{remark}
\newtheorem{remark}[theorem]{Remark}
}

{\theoremstyle{remark}

}

\numberwithin{equation}{section}

\renewcommand{\labelenumi}{{ \theenumi.}}
\renewcommand{\labelenumii}{{(\alph{enumii})}}

\def\la{\lambda}
\def\al{\alpha}

\def\choose #1 #2{\begin{pmatrix}#1\\#2\end{pmatrix}}

\def\s{{\bf s}}

\begin{document}

\title{On representations of $U'_q\so_n$}

\author{Hans Wenzl}

\address{Department of Mathematics\\ University of California\\ San Diego,
California}

\email{hwenzl@ucsd.edu}

\begin{abstract} We study representations of the non-standard quantum deformation $U'_q\so_n$ of $U\so_n$ via a Verma module approach.
This is used to recover the classification of finite-dimensional modules
for $q$ not a root of unity, given by classical and non-classical series.
We obtain new results at roots of unity, including the classification
of self-adjoint representations on Hilbert spaces.
\end{abstract}

\maketitle

\vskip .5cm

It is well-known that there exists another $q$-deformation $U'_q\so_n$
of the universal enveloping algebra of $\so_n$ which differs from
the Drinfeld-Jimbo quantum group $U_q\so_n$, see \cite{GK}, \cite{L}, \cite{NS}.
Unlike the latter one, the algebra $U'_q\so_n$ can be
embedded as a coideal subalgebra into the Drinfeld-Jimbo quantum group
 $U_q\sl_n$. Besides its occurrences
in connection with quantum symmetric spaces and in mathematical
physics, it has more recently  appeared as a centralizer algebra
for tensor products of spinor representations of $U_q\so_N$ (\cite{WSp})
and for $q$-Howe duality for orthogonal quantum groups (\cite{ST}),
and also in connection with von Neumann subfactors and in quantum computing
(\cite{WFu}, \cite{RW}). In
the last two areas mentioned, one has to deal with representations
of $U'_q\so_n$ at roots of unity. A detailed classification of
simple finite-dimensional representations $U'_q\so_n$
has been obtained in the
work of Klimyk and his collaborators if $q$ is not a root of unity
(see \cite{IK} and references therein). Some representations have
also been obtained for $q$ a root of unity in \cite{IK2}.
Unfortunately, these
representations are not the ones which appear
in connection with subfactors and quantum computing.
This lead to the approach in this paper, which differs from the ones
by other authors in several substantial ways
(see Section \ref{knownrep} for previous approaches, and the next paragraph for an outline of the
approach in this paper as well as the discussion in the last section).
It is also hoped that at least parts of
this approach may be useful for the study
of representations of other coideal algebras, see e.g.
\cite{ES}, \cite{ST}, \cite{LC} and the discussion at the
end of this paper.

The first step in our approach is to define and construct  Verma
modules.
Unfortunately, the algebras $U'_q\so_n$ do not have $canonical$
raising and lowering operators, which makes the construction more
difficult (see Section \ref{discussresults} for details).
Nevertheless, one can
construct an analog of a Verma module $V_m$ also for the algebra
$U'_q\so_n$, for arbitrary weight $m$. As a first application,
we obtain the classification of finite-dimensional simple modules
for $q$ not a root of unity. Unlike for
Drinfeld-Jimbo quantum groups, there also exist finite dimensional
simple representations of $U'_q\so_n$ which are not deformations of
representations of $U\so_n$; following the notation in \cite{IK}, we
call them non-classical representations.
Moreover, we construct canonical finite-dimensional  quotients
in these cases which are well-defined for all values of $q\neq 0$.
Their dimensions are given by Weyl's character formula. If $q$ is
not a root of unity,
they are irreducible for the classical modules, but reducible for
the non-classical modules. In the latter case, the quotient decomposes
into the direct sum of $2^{\lfloor (n-1)/2\rfloor}$ mutually non-equivalent
irreducible modules,
all of which have the same character with respect to our chosen
Cartan subalgebra.
With this approach, one can now also extend our classification
to certain classes of representations at roots of unity,
including all representations
on Hilbert spaces on which the generators act as self-adjoint operators.
These representations are different from the ones constructed in \cite{IK2}.

Here is the contents of our paper in more detail. After basic
definitions, we first give a detailed discussion of known results about
representations of $U'_q\so_n$ in Section \ref{knownrep}.
In the second
section we define the notion of a Verma module for $U'_q\so_n$, and
we determine a certain canonical spanning set. In the third
section, we prove all the necessary results for Verma modules for
$U'_q\so_3$ via elementary methods which more or less have been
known before. In particular, we recuperate the classification of
all finite-dimensional $U'_q\so_3$-modules.
Using this and known results about representations of $U'_q\so_n$,
as reviewed in the first section, we prove that the spanning set
in the second section is indeed a basis for any Verma module.
Section 5 starts out
with an elementary study of representations of $U'_q\so_4$. This is
then used to classify all weights for which the corresponding Verma
module allows a finite-dimensional  quotient, and to construct such
a quotient. These modules can be viewed as analogs of (dual) Weyl
modules and their dimensions are given by Weyl's character formula.
As already mentioned above, they are reducible in the case of
non-classical representations.
In Section 6, we apply our results for representations at
roots of unity. In particular, we show that for representations on
Hilbert spaces for which the generators of $U'_q\so_n$ act as
self-adjoint operators, they are again classified by their highest
weights. The final section contains a brief discussion comparing the
approach in this paper with other approaches and possible applications.

{\it Acknowledgements} I would like to thank the referees for their careful
reading of  previous versions of this manuscript which helped to substantially improve
the presentation of the material. Part of the work on this manuscript was
done while I enjoyed the hospitality of the Center for Quantum Geometry
of Moduli Spaces in Aarhus.

\section{Definitions and known representations}\label{defrep}

\subsection{Definitions}\label{definitions} We are primarily interested
in representations over the complex numbers, with $q$ also a complex number.
Occasionally, it will be convenient to
view $q$ as a variable over the complex numbers.
E.g. the results in Section \ref{spanset} can be formulated over the ring of Laurent polynomials
in $q$.
The algebra $U_q'\so_n$, often referred to as non-standard
deformation of the universal enveloping algebra $U\so_n$ of the orthogonal Lie algebra $\so_n$ (see \cite{GK}, \cite{L}, \cite{NS}),
is defined via generators $B_i$, $1\leq i\leq n-1$ and relations $B_iB_j=B_jB_i$ for $|i-j|>1$ and
\begin{equation}\label{basicrel}
B_i^2B_{i\pm 1}-[2]B_iB_{i\pm 1}B_i + B_{i\pm 1}B_i^2=B_{i\pm 1},
\end{equation}
where the same sign is chosen in each summand.
We identify roots and weights of $U\so_n$ with
vectors in $\R^k$, where $k=n/2$ or $(n-1)/2$ depending on the parity of $n$,
as usual. We shall use the notation $\ep_i$
for the $i$-th standard basis vector of $\R^k$.
As usual, we denote the simple roots $\al_i$ for $\so_n$ by
$\al_i=\ep_i-\ep_{i+1}$ for $1\leq i\leq (n-2)/2$, and by
$\al_k=\ep_{k-1}+\ep_k$ for $n=2k$ even, and $\al_k=\ep_k$ for $n=2k+1$ odd.

The analog of the Cartan subalgebra in $U_q'\so_n$ is the algebra generated by
$B_1, B_3, ... B_{2k-1}$ for $n=2k$ or $n=2k+1$. A vector $v$ in a $U_q'\so_n$-module is said to have weight
$m$ if $B_{2i-1}v=m_iv$
for all $B_{2i-1}\in \h$. If $\la\in\R^k$ we shall often use the notation
$[\la]$ for the weight whose coordinates are $[\la_i]$;
here, as usual, $[n]=(q^n-q^{-n})/(q-q^{-1})$. Another important type of
weights are denoted by $[\la]_+$, where the $i$-th coordinate is given
by $[\la_i]_+$, where now $[n]_+=\sqrt{-1}(q^n+q^{-n})/(q-q^{-1})$.

\subsection{Known representations}\label{knownrep}
We review some of the known representations of $U'_q\so_n$. At this point, it mainly serves as a motivation for this paper. The results of this
subsection will only be used later when we prove the linear independence of the spanning set of a Verma module, as established in the next section.

1. A classification of finite dimensional representations of $U'_q\so_n$
for $q$ not a root of unity has been given by Klimyk and his
coauthors, see \cite{IK} and the literature quoted there.
In particular, they found, besides the expected deformations of
representations of $U\so_n$ (called classical representations),
additional representations which are not well-defined
in the classical limit $q=1$. Not surprisingly, the classical representations
are classified by the dominant integral weights of $\so_n$.
We shall reprove
their results (or results equivalent to theirs) in the special cases of
$U'_q\so_3$ and $U'_q\so_4$, essentially in the same form.
In the approach of \cite{IK} the authors then give explicit matrix
representation for a Gelfand-Tseitlin type basis for representations in the
general case, using tensor operators and a $q$-version of the Wigner-Eckart theorem.
The matrix coefficients in these explicit
representations are considerably more complicated in the higher rank cases.
Moreover, these representations are not well-defined for $q$ a root of unity if the highest weight $\la$ is
an integral dominant weight. We will define
representations for $n>4$ in a different way for both the classical
and the non-classical representations.

2. Quite different representations have been found in \cite{WSp}. Let $S$
be the spinor representation of the Drinfeld-Jimbo quantum group $U_q\so_N$.
Then a matrix $C\in \End(S^{\otimes 2})$ was found in \cite{WSp}
such that the map
\begin{equation}\label{WSprep}
B_i\ \mapsto \ 1_{i-1}\otimes C\otimes 1_{n-1-i}\ \in\ \End(S^{\otimes n})\quad 1\leq i\leq n-1
\end{equation}
defines a representation of $U'_q\so_n$. In particular, any irreducible
classical representation of $U'_q\so_n$ does appear in these representations,
see e.g. \cite{RW}, Theorem 2.2. Moreover, the relations are much easier to check as for the
representations discussed before, and the representations are also
well-defined for roots of unity. In particular, it follows from \cite{Wcat}
and \cite{WSp} that for special roots of unity
(usually of the form $q=e^{\pm 2\pi i/\ell}$)
one can obtain representations of $U'_q\so_n$ on Hilbert spaces
on which the generators $B_i$ act via self-adjoint operators.
We may sometimes call such representations of $U'_q\so_n$ $unitary$, as they
are (at least for these specific examples) closely related to unitary
representations of braid groups. 

Representations of similar type have also been found independently by Sartori and
Tubbenhauer in \cite{ST} by a different method. Instead of $S$, they consider
a $q$-version of the exterior algebra of the vector representation of $\so_N$.
It is isomorphic to $S^{\otimes 2}$ as a $U_q\so_N$-module for $N$ even.
In this case, their result is more or less equivalent to the one in \cite{WSp}
for even tensor powers of $S$.  

As non-classical representations are not well-defined for $q=1$, we do not
expect them to appear in the papers just mentioned. However, after having
submitted the first version of this paper, the author found that one can also
obtain non-classical representations of $U'_{\tilde q}\so_n$ as 
commutants of $U_q\so_N$ on tensor powers of  spinor representations
for a suitably modified $\tilde q\neq q$. Unlike for classical representations,
the parameter $\tilde q$ for the representations of $U'_{\tilde q}\so_n$
will go to $-1$ for the classical limit $q\to 1$. In particular, one also
obtains unitary non-classical representations of $U'_{\tilde q}\so_n$ for roots of unity $\tilde q$ near
$-1$.  This is currently being worked out
in \cite{Wdual}.

3. Representations of $U'_q\so_n$ for $q$ a primitive $\ell$-th root of unity
also appear in \cite{RW}, Section 4. We define them for $n$ odd,
with the even case obtained by restricting to $U'_q\so_{n-1}$.
Let $V$ be an
$\ell^{(n-1)/2}$-dimensional vector space with basis $v(\vec i)$, where
$\vec i\in \{0,1,\ldots, \ell-1\}^{(n-1)/2}$.
The action of $u_{2s-1}$
on $V$ is defined by $u_{2s-1}v(\vec i)=q^{i_s}v(\vec i)$.
The action of $u_{2s}$ is given  by the rule (indices modulo $\ell$):

$$u_{2s}(v(i_1,\ldots,i_{s},i_{s+1},\ldots,i_{\frac{n-1}{2}}))=v(i_1,\ldots,
i_s+1 , i_{s+1}-1,i_{s+2},\ldots,i_{\frac{n-1}{2}});$$
in other words, the even indexed generators $u_{2s}$ permute the vectors
$v(i_1,\ldots,i_{\frac{n-1}{2}})$
by shifting the $s$-th index up by 1 and the
$(s+1)$-st index down by 1, except for $s=(n-1)/2$ where there is no  index
left for shifting down. It is easy to check that these operators satisfy the relations
 $u_iu_{i+1}=qu_{i+1}u_i$ and $u_iu_j=u_ju_i$ for $|i-j|>1$.
It was shown in \cite{RW}, Lemma 4.3 that we obtain two representations
of $U'_q\so_n$ on $V$ given by the maps

\begin{equation}\label{qtorus1}
B_j\mapsto b_j=\pm \frac{\sqrt{-1}}{q-q^{-1}}\ (u_j+u_j^{-1})
\quad
1\leq j\leq n-1,
\end{equation}
\begin{equation}\label{qtorus2}
B_j\mapsto b_j= \frac{u_j+u_j^{-1}}{q-q^{-1}},\quad
1\leq j\leq n-1.
\end{equation}

\begin{remark}\label{motivation} The main motivation for this paper came
from the paper \cite{RW}. One can obtain
 unitary braid representations from the representations mentioned under 2 as well
as under 3 above which are of interest in quantum
computing. It was conjectured that the braid representations from Example 3
correspond to certain cases within the framework of Example 2.
The proof of showing this correspondence was reduced to showing that
the corresponding representations of $U'_q\so_n$ are isomorphic.
However, as we could not find a general representation theory which would
also involve roots of unity, this was finally done via a somewhat indirect
procedure. One of the aims of this paper is to close this gap by extending
existing results also to roots of unity, at least to the extent that it
would cover the just mentioned examples.
\end{remark}

\section{Spanning set for Verma modules}\label{spanset}

\subsection{Basic exchange relations} It has already been observed in
\cite{GK} that a PBW type theorem holds for the algebra $U'_q\so_n$, using
its embedding into the quantum group $U_q\sl_n$. We will give a direct proof
here of a weaker version of this statement, by exhibiting an explicit spanning set.
This will then be used to construct an explicit spanning set for
an analog of a Verma module. We view $q$ as a variable in this section.

\begin{lemma}\label{lemrelations} We define $B_{n,k}=B_{n-1}B_{n-2}\ ...\ B_k$.
Then we have the following equalities, where $w_1$, $w_2$ and $w_3$ are linear combinations of products of generators with coefficients being Laurent polynomials in $q$:

(a) \ $B_{n,k}B_k=[2]B_kB_{n,k}-B_k^2B_{n,k+1}+B_{n,k+1}$, if $n>k+1$,

(a)$'$\ $B_{n-1}B_{n,k}=[2]B_{n,k}B_{n-1}-B_{n-1,k}B_{n-1}^2+B_{n-1,k}$, if $k<n-1$,

(b) \ if $k<r<n-1$, we have
$$[B_r,B_{n,k}]=B_{n,r+2}[B_r, B_{r,k}B_rB_{r+1}]=[B_r,B_{r,k}B_rB_{n,r+1}]=[B_{r-1}B_rB_{n,r+1},B_r]B_{r-1,k},$$

(c) \ $B_rB_{r+1}B_{n,k}= w_1B_{r+2}+w_2B_kB_{r+1}+w_3B_r$.
\end{lemma}

$Proof.$ These are mostly straightforward calculations, where (a) and $(a)'$ follow fairly easily from the generating relations.
The proof of (b) uses the following calculation, where we use $(a)$ and $(a)'$:
\begin{align}
B_{n,r-1}B_r\ & =\ B_{n,r+1}B_rB_{r-1}B_r\ =\ \frac{1}{[2]}B_{n,r+1}(B_r^2B_{r-1}+B_{r-1}B_r^2-B_{r-1})\ =\ \cr
& = \ \frac{1}{[2]}([2]B_rB_{n,r}B_{r-1}-B_r^2B_{r-1}B_{n,r+1}+B_{r-1}B_{n,r}B_r)\ =\cr
& =\ B_rB_{n,r}B_{r-1}+ B_{r-1}B_rB_{n,r}-\frac{1}{[2]}(B_r^2B_{r-1}+B_{r-1}B_r^2-B_{r-1})B_{n,r+1}\cr
& =\ B_rB_{n,r-1}-[B_r,B_{r-1}B_rB_{n,r+1}].
\end{align}
A similar calculation also shows that $$[B_{n-2},B_{n,k}]=[B_{n-2},B_{n-2,k}B_{n-2}B_{n-1}].$$
The general case for statement (b) can be reduced to these calculations by observing that  $B_rB_{n,k}=
B_{n,r+2}(B_rB_{r+2,k})$ and $B_{n,k}B_r=(B_{n,r-1}B_r)B_{r-2,k}$.
Statement (c) can now be deduced from (a) and (b). $\Box$

\subsection{Ordering}
We will need an ordering on monomials in the generators which is defined as follows: If two monomials $w_1$ and $w_2$ have different
lengths, i.e. the number of factors is different, then the word of shorter length is the smaller one. For words of same length,
it will be convenient to
define inverse alphabetic order, reading from the right. E.g. we get
for words of length 3 in $U'_q\so_3$ the ordering
$$B_2^3<B_1B_2^2<B_2B_1B_2<B_1^2B_2<B_2^2B_1<B_1B_2B_1<B_2B_1^2<B_1^3.$$
Informally, the more generators with high indices are to the right in a monomial,
the smaller it is.
\begin{lemma}\label{wordord} Let $k>l$. Then $B_{n,k}B_{n,l}=[2] B_{n,l}B_{n,k} + $ lower terms.
\end{lemma}

$Proof.$ The proof goes by downward induction on $k$. If $k=n-1$, then $B_{n,n-1}=B_{n-1}$, and the claim follows from
Lemma \ref{lemrelations}, (a)$'$. For the induction step, we use Lemma \ref{lemrelations}, (b) for the second line below,
and the induction assumption for the third line.
\begin{align}
B_{n,k}B_{n,l}\ &=\ B_{n,k+1}B_kB_{n,l}\ =\ \cr
&=\ B_{n,k+1}B_{n,l}B_k + B_{n,k+1}(B_{n,k+2}[B_k,B_{k,l}B_kB_{k+1}])\ =\cr
&=\ ([2] B_{n,l}B_{n,k+1}+ {lower\ terms})B_k+ B_{n,l}B_kB_{n,k+1}-B_{k,l}B_{n,k}B_{n,k}.
\end{align}
It follows directly from the definitions for words $w_1$ and $w_2$ in the generators $B_j$ that if $w_1<w_2$,
then also $w_1B_k<w_2B_k$. The claim follows from this. $\Box$
\vskip .3cm

\subsection{A spanning set for a Verma module}
Let $m$, $\tilde n$ be weights. Then we define the left ideal $I_{m,\tilde n}$ by
\begin{equation}\label{Ilambdadef}
I_{m,\tilde n}=U'_q\so_n\langle (B_{2i-1}-m_i1), (B_{2i-1}B_{2i}-\tilde n_iB_{2i})\rangle
\end{equation}
for all values of $i$ for which the indices $2i-1$ and $2i$ are between
(including) 1 and $(n-1)$. We will show later, see Lemmas \ref{basiclem} and  \ref{weightshift}, 
that for given $m_i$ we obtain a nonzero vector $B_{2i-1}B_{2i}$ mod $I_{m,\tilde n}$
only for two special values of $\tilde n_i$. In the classical case these two cases would correspond to 1 being a highest respectively 
a lowest weight vector mod $I_{m,\tilde n}$. This will be made more precise in the following sections.

\begin{proposition}\label{Vermaprop} Let $R$ be the ring of Laurent polynomials in $q$,
and view $U'_q\so_n$ as an $R$-algebra. Moreover, let $R[m,\tilde n]$ be the polynomial ring over the 
coordinates $m_i$ and $\tilde n_i$ of the weights $m$ and $\tilde n$.

(a) The algebra $U'_q\so_n$ is spanned over $R$ by monomials of the form
$$B_{2,1}^{e(2,1)}B_{3,1}^{e(3,1)}B_{3,2}^{e(3,2)}B_{4,1}^{e(4,1)}\ ...\ B_{n,n-2}^{e(n,n-2)}B_{n,n-1}^{e(n,n-1)},$$
where the $e(i,j)$ are nonnegative integers.

(b) The quotient $V_{m,\tilde n}=U'_q\so_n/I_{m,\tilde n}$
is spanned over  $R[m,\tilde n]$ by the set 
$$\B\ =\ \{ B_{3,2}^{e(3,2)}B_{4,2}^{e(4,2)}B_{5,2}^{e(5,2)}B_{5,4}^{e(5,4)}\ ...\ B_{n, 2\lfloor (n-3)/2\rfloor}^{e(n, 2\lfloor (n-3)/2\rfloor)}B_{n, 2\lfloor (n-1)/2\rfloor}^{e(n, 2\lfloor (n-1)/2\rfloor)}+I_{m,\tilde n}\},$$
i.e. by the residue classes mod $I_{m,\tilde n}$ of elements  in (a) which only contain factors $B_{k,r}$ for which
$r$ is even.
\end{proposition}

$Proof.$ We show (a) by induction on $n$. It is clear that $U'_q\so_2$ is spanned by the powers of $B_1$.
For the induction step $n$ to $n+1$ it suffices to show that $U'_q\so_{n+1}$
is spanned over $R$ by elements of the form
$$x\ B_{n+1,1}^{e(n+1,1)}B_{n+1,2}^{e(n+1,2)}\ ...\ B_{n+1,n-1}^{e(n+1,n-1)}B_{n+1,n}^{e(n+1,n)},$$
with $x\in U'_q\so_n$.  In the first step, assume we have a monomial containing a string $B_{n+1,k}B_r$,
with $k\leq r<n$. But then we can express this string by a linear combination over $R$ of smaller monomials
using Lemma \ref{lemrelations}, (a) if $r=k$, and the last expression in the statement of  Lemma \ref{lemrelations}, (b)
for $r>k$. Indeed, it is easy to see that the generator $B_n$ has moved to the right in all of these expressions.
Hence we can assume that elements of the form $x\ \prod_{i=1}^b B_{n+1,k_i}$, $x\in  U'_q\so_n$
span $U'_q\so_{n+1}$.
To complete the proof, it suffices to show that we can replace an expression of the form $B_{n+1,k_1}B_{n+1,k_2}$
with $k_1>k_2$ by a linear combination over $R$ of smaller monomials. But this follows from Lemma \ref{wordord}.

To prove part (b), we will show that whenever a product $\prod_i B_{n_i,k_i}$
contains an odd $k_i$, we can replace it by a linear combination of products
of generators  modulo $I_{m,\tilde n}$ where each product has fewer factors than the original element.
For the purpose of the induction, we will prove such a statement by induction on $s$ for
any element of the form $B_kw$ or $B_kB_{k+1}w$, with $k$ odd and
$w=\prod_{j=1}^sB_{n_j,k_j}$ with all $k_j$ even.
If $s=0$, this follows directly from the definition of $I_{m,\tilde n}$.
For the induction step, we only need to observe that we can replace
$B_kB_{n_1,k_1}$ respectively $B_kB_{k+1}B_{n_1,k_1}$
by suitable linear combinations of
elements ending with an element $C$, where $C$ is equal to $B_k$, $B_kB_{k+1}$ or $B_{k+2}$, by Lemma \ref{lemrelations}.
Now, by induction assumption, we can replace $C\prod_{j=2}^sB_{n_j,k_j}$ by a linear combination of shorter elements. $\Box$

\medskip
Let $\s=(s_i)$, with $s_i$ non-negative integers for $1\leq i\leq \lfloor (n-1)/2\rfloor$. We define the subspace $V_{m,\tilde n}[[\s ]]$ to
be the span of all monomials in $\B$ which contain the generator $B_{2i}$
at most $s_i$ times. The proof of the following corollary is
essentially just a special case of the proof of Proposition \ref{Vermaprop},(b).

\begin{corollary}\label{Vermapropcor}
The subspace $V_{m,\tilde n}[[\s ]]$ is a module of the Cartan subalgebra of $ U'_q\so_n$,
i.e. the subalgebra generated by the elements $B_{2i-1}$, $i\leq n/2$.
\end{corollary}

\section{Representations of $U'_q\so_3$}\label{so3repsec}

As in the classical case, the representation theory for the simplest
nontrivial case, $U'_q\so_3$ is both elementary and important. Most
or probably all of the results in this section have been obtained
before, see e.g. \cite{IK} and references there, or \cite{HP}. We
reprove the results here, as it will fix our notations and the
applied methods will be useful for the general case.

\subsection{General weights and Verma modules}
As usual, we identify weights of $U'_q\so_3$ with eigenvalues of $B_1$.
So the ideal $I_{m,n}$ defined in the last section would depend on two
numbers $m$ and $n$. We also use the
notation $v_0$ for the image of 1 in the quotient $V_{m,n}=U'_q\so_3/I_{m,n}$,
and the notation
\begin{equation}\label{weightspace}
V_{m,n}[k]={\rm span} \{ B_2^jv_0,\ 0\leq j\leq k\ \}.
\end{equation}
In order to describe the weights of $V_{m,n}$, we define, for given weight $m$,
a sequence $(m_j)$ by $m_0=m$,
\begin{equation}\label{defmj}
m_{1}=\frac{1}{2}([2]m\pm \sqrt{(q-q^{-1})^2m^2+4}), \quad {\rm and}\quad m_{j+1}=[2]m_j-m_{j-1}
 \quad {\rm for\ }j>0.
\end{equation}
This sequence is uniquely determined by $m$ and $m_1$, where $m_1$ is one of two roots of the polynomial $x^2-[2]mx+m^2-1$. 
Such sequences will only appear in this section. So no confusion with coordinates $m_i$ of a weight $m$
for higher rank cases should occur.
The following examples of sequences $(m_j)$ will be particularly important
in our paper:

(a) If $m=[\la]=\frac{q^\la-q^{-\la}}{q-q^{-1}}$, then
$m_j=[\la\pm j]$, with the sign fixed by the choice of sign for $m_1$.

(b) If $m=[\la]_+=i\frac{q^\la+q^{-\la}}{q-q^{-1}}$, then
$m_j=[\la\pm j]_+$, with the sign fixed by the choice of sign for $m_1$.

\begin{lemma}\label{mbasic} (a) The numbers $m_{j\pm 1}$ are the zeros of $x^2-[2]m_jx+m_j^2-1$ for all $j>0$.
In particular, $m_{j-1}\neq m_{j+1}$ unless $m_j=\pm [0]_+=\pm 2i/(q-q^{-1})$.

(b) We have $m_{j+k}=m_j$ for some $k>0$ only if $q^k=1$ or $m_j=\pm [k/2]_+$
\end{lemma}

$Proof.$ Part (a) is shown by induction on $j$, with $j=0$ being trivially
true. It follows from the induction assumption for $j$ that
$m_{j+1}^2-[2]m_jm_{j+1}+m_j^2-1=0$; it follows
that $m_j$ is a root of the
polynomial $x^2-[2]m_{j+1}x+m_{j+1}^2-1$.
But then the second root is equal to $[2]m_{j+1}-m_j=m_{j+2}$.
The second statement follows from setting the discriminant of this quadratic equation equal to 0.
To prove (b), let us assume $j=0$ for ease of notation. It is easy
to prove by induction on $k$, using $[2]m_j=m_{j+1}+m_{j-1}$ that
\begin{equation}\label{minduction}
m_k=[k]m_1-[k-1]m_0\quad {\rm and} \quad m_{-1}=[k+2]m_k-[k+1]m_{k+1}.
\end{equation}
Solving for $m_1$ in the first formula, we obtain from $m_k=m_0$ that
$$(2+2[k-1]-[2][k])m_0=\pm [k]\sqrt{(q-q^{-1})^2m_0^2+4}.$$
This can be transformed to
$$(2-q^k-q^{-k})m_0^2=[k]^2.$$
Hence, if $0\neq -(q^{k/2}-q^{-k/2})^2=-q^{-k}(1-q^k)^2$, then $m_0=\pm [k/2]_+$. $\Box$

\begin{lemma}\label{basiclem} (a) If $v$ is a vector in the
$U'_q\so_3$ module $V$ of weight $\mu$,
then $(B_1^2-[2]\mu B_1+(\mu^2-1)B_2v=0$. In particular, if $\mu\neq \pm [0]_+$,
$B_2v$ is a linear combination of two weight vectors.

(b) The quotient $U'_q\so_3/I_{m,n}$ can have dimension $>1$
only if $n$ is a root of the
polynomial $x^2-[2]mx+ m^2-1$.

(c) $B_1B_2^kv_0=m_kB_2^kv_0+v'_k$, where $v'_k\in V_{m,n}[k-1]$,
as defined in \ref{weightspace}

(d) We have $\prod_{j=0}^k (B_1-m_j)=0$ on $V_{m,n}[k]$.

\end{lemma}

$Proof.$ Claim (a) follows from
$B_1^2B_2v=([2]B_1B_2B_1-B_2B_1^2+B_2)v$, using $B_1v=\mu v$,
and from Lemma \ref{mbasic}, (a). As
$B_2$ mod $I_{m,n}$ would be an eigenvector of $B_1$ with eigenvalue
$n$, it can only be nonzero mod $I_{m,n}$ if $n$ is a root of the
polynomial as stated, by (a). If $n$ is not a root of the polynomial,
it follows that $B_2^j=0$ mod
$I_{m,n}$ for all $j\geq 1$. Hence the quotient has at most
dimension 1. Claim (c) is obviously true for $k=0,1$. For $k>1$, it
follows by induction, using
$$B_1B_2^k=[2]B_2B_1B_2^{k-1}-B_2^2B_1B_2^{k-2}+B_1B_2^{k-2}.$$
Assuming that the vectors $B_2^j$, $0\leq j\leq k$ are linearly
independent,
it follows from part (c) that $B_1$ acts on $V_{m,n}[k]$
via a triangular matrix with diagonal entries $m_j$, $0\leq j\leq k$.
If $B_2^{k+1}v_0\subset V_{m,n}[k]$, then so are all higher powers of $B_2$,
and the same argument can be used for the smallest $k$ for which
the $B_2^jv_0$ are linearly independent, $0\leq j\leq k$. $\Box$

\subsection{Finite-dimensional modules} The statement of the
following proposition is not true for $q\neq 1$ a root of unity.

\begin{proposition}\label{finVerm}
If $q$ is not a root of unity, then any finite dimensional simple $U'_q\so_3$-module has to be a quotient of a Verma module.
\end{proposition}

$Proof.$ Let $V$ be a finite dimensional simple $U'_q\so_3$-module. Assume we can find
an eigenvector $v$ of $B_1$ with eigenvalue $m_0$ such that also $B_2v$ is an eigenvector of $B_1$.
Then the eigenvalue of $B_2v$ would have to be one of the two possible values $m_1$ given in
\ref{defmj}. It follows that the map $u\in U'_q\so_3\mapsto uv$ has a kernel containing the ideal $I_{m_0,m_1}$. 
As $V$ is simple, this map would be surjective, and the claim would follow.

By Lemma \ref{basiclem}, (a), the vector $B_2v$ lies in an at most
2-dimensional $B_1$-invariant subspace for any weight vector $v$. Let us first assume that our module 
does not contain the weight $\pm [0]_+$. If we can not find a $v$ as
in the previous paragraph, we can inductively construct a sequence
of eigenvectors $(v_j)$, $j\in\Z$ with eigenvalues $m_j$, by Lemma
\ref{basiclem}, (a) and Lemma \ref{mbasic}, (a). 
As $V$ is finite
dimensional, this sequence can only take finitely many values and we
have $m_j=m_{j+k}$ for some $j$ and some $k>0$. If $q$ is not a root
of unity, then $m_j=\pm [k/2]_+$, by Lemma \ref{mbasic}, (b). But then
$m_s=\pm [j+k/2-s]_+$, by Lemma \ref{mbasic}, (a) and statement (b) before it,
and the sequence $(m_s)$ would take infinitely many
distinct values. This leads to a contradiction to $V$ being finite dimensional.
A similar argument also works if $\pm [0]_+$ is a weight of $V$, or see Lemma \ref{mplus}. $\Box$

\subsection{Basis for Verma module}  We have seen in Lemma \ref{basiclem}
that, for given $m$, we only obtain
nontrivial quotients for two values of $n$. As a consequence of this,
we will  simplify the notation by writing $I_m$ for $I_{m,n}$,
and $V_m=U'_q\so_3/I_m$, $V_m[k]$ etc,
with the understanding that one of the two
possible choices $m_1$ for $n$ has been fixed.

\begin{proposition}\label{weightbasis} (a) Assume $m=m_0$ is such that
$m_j\neq\pm [0]_+=\pm\frac{2i}{q-q^{-1}}$ for all $j>0$, and that $q$ is not a root of unity. Then the Verma
module $V_m$ has a basis of weight vectors $v_j\in V_m[j]$ with
weight $m_j$, $j\geq 0$. We can choose them such that
\begin{equation}\label{Baction}
B_2v_j=v_{j+1}+\al_{j-1,j}v_{j-1},
\end{equation}
where
$$\al_{j-1,j}=\frac{m_0m_{-1}-m_{j-1}m_j}{(m_{j-1}-m_{j+1})(m_{j-2}-m_j)}.$$
(b) We are not in case (a) if and only if $m=\pm [\la]_+$ for a
positive integer $\la$ and $m_1=\pm [\la-1]_+$. In this case, there
still exist weight vectors $v_k$ for $0\leq k\leq \la$ and for $k\geq 2\la +1$,
with $m_j=\pm [\la-j]_+$. Here the signs equal the one for $m$.
\end{proposition}

$Proof.$ Let us first assume that the vectors $B_2^jv_0$, $j=0,1,2,\ ...$
are linearly independent. Then it follows from
Lemma \ref{basiclem}, (c) and (d) that $B_1$ acts via a triangular matrix
with respect to this basis, with diagonal entries $m_j$.
It follows from Lemma \ref{basiclem}, (a) by induction on $j$ that
$B_1$ can be diagonalized, as $m_{j+1}\neq m_{j-1}$ for all $j>0$.
Observe that this also works if $m_0=\pm [0]_+$.

We normalize the eigenvector $v_k$ with eigenvalue $m_k$
such that $v_k=B_2^kv+\sum_{j=0}^{k-1} \beta_jB_2^jv$. From this follows
the expression \ref{Baction} for $B_2v_j$, up to determining the scalar $\al_{j-1,j}$.
To do this, first observe that
\begin{equation}\label{prop32}
B_2^2v_j=v_{j+2}+(\al_{j,j+1}+\al_{j-1,j})v_j +\al_{j-2,j-1}\al_{j-1,j}v_{j-2}.
\end{equation}
We now calculate $B_1B_2^2v_j$ in two different ways. First by applying $B_1$
to  Eq \ref{prop32} directly, using $B_1v_i=m_iv_i$.
Secondly, using the relations, we also obtain
$$B_1B_2^2v_j= [2]B_2B_1B_2v_j-B_2^2B_1v_j+B_1v_j,$$
after which we again expand it into a linear combination of $v_{j+2}$, $v_j$ and $v_{j-2}$
using the formula \ref{Baction} for $B_2v_j$.
Comparing the coefficients of $v_j$ in these two expressions we obtain
$$m_j(\al_{j,j+1}+\al_{j-1,j})= [2](m_{j+1}\al_{j,j+1}+m_{j-1}\al_{j-1,j})-m_j(\al_{j,j+1}+\al_{j-1,j}) +m_j.$$
Using $m_{j+2}=[2]m_{j+1}-m_j$, we can simplify the expression above to
\begin{equation}\label{recursealpha}
(m_{j+2}-m_j)\al_{j,j+1}=(m_j-m_{j-2})\al_{j-1,j}-m_j.
\end{equation}
Hence we can express $\al_{j,j+1}$ in terms of $\al_{j-1,j}$. To get the induction going,
we calculate $B_1B_2^2v_0$ in two different ways as before. One checks directly that
$$\al_{0,1}=\frac{m_0}{m_0-m_2}=\frac{m_0(m_{-1}-m_1)}{(m_0-m_2)(m_{-1}-m_1)}.$$
The general form for $\al_{j,j+1}$ now follows from the two previous formulas by
induction on $j$.

To prove the linear independence of the vectors $B_2^jv_0$, we consider
a vector space $V$ with basis $v'_j$, on which we define the action
of $B_1$ and $B_2$ by $B_2v_j'=v'_{j+1}+\al_{j-1,j}v'_{j-1}$ and by $B_1v'_j=m_jv_j'$. It follows
essentially from the same calculations as before that this does define
a representation of $U'_q\so_3$. It is now straightforward to check
that the kernel of the map $u\in U'_q\so_3\mapsto uv'_0$ contains the ideal
$I_m$, and that the vectors $B^j_2v'_0$ are linearly independent.

Finally, to check (b), assume $m_j=\frac{2i}{q-q^{-1}}=[0]_+$. One deduces from
this by induction on $r$ that $m_{j-r}=[r]_+=[-r]_+$, as $[n]_+=[-n]_+$, 
and hence $m=[j]_+$.
For the last claim, one shows the existence of eigenvectors $v_k$ for $k\leq j$
as in case (a).
One observes that the matrix representing $B_1$ on $V_m[2j+1]$ for $m=[j]_+$
has the eigenvalue $[-j-1]_+=[j+1]_+$ with multiplicity 1, which does not exist for its
restriction to $V_m[ 2j]$, by Lemma \ref{basiclem}, (d).
Hence the projection onto the eigenspace of $B_1$ corresponding
to eigenvalue $[\la-2r-1]_+$ is well-defined for $\la=r$, and
$v_{2r+1}$ is its unique eigenvector with coefficient 1 for
$B_2^{2r+1}v_0$ in its expansion with respect to the basis
$(B_2^jv_0$, $0\leq j\leq 2r+1)$. The existence of eigenvectors $v_k$ with $k>2j+1$
 is again shown as in (a), using Lemma \ref{basiclem}. $\Box$

\begin{corollary}\label{ali} (a) If $m=[\la]$ and $m_1=[\la-1]$, then
$$\al_{j-1,j}=\frac{[2\la+1-j][j]}{(q^{\la-j}+q^{j-\la})(q^{\la-j+1}+q^{j-1-\la})}.$$

(b) If $m=[\la]_+$ and $m_1=[\la-1]_+$, then
$$\al_{j-1,j}=\frac{[2\la+1-j][j]}{(q^{\la-j}-q^{j-\la})(q^{\la-j+1}-q^{j-1-\la})}.$$
\end{corollary}

\begin{corollary}\label{so3pol} Define polynomials $P_n$ inductively by
$P_0=1$, $P_1(x)=x$ and $P_{n+1}(x)=xP_n(x)-\al_{n-1,n}P_{n-1}(x)$.
Then $P_n(B_2)v_0=v_n$ for all $n\geq 0$.
\end{corollary}

\subsection{Highest weight vectors} The following lemma will be useful
for determining non-trivial highest weight vectors in $V_m$.

\begin{lemma}\label{finitelem} Assume $\al_{j,j+1}=0$.

(a) If $j=2k$ is even, then $m_k=0$ or $m_k=\pm [0]_+=\pm 2i/(q-q^{-1})$.

(b) If $j=2k+1$ is odd, then $m_k=\pm [1/2] $ or $m_k=\pm [1/2]_+$.
\end{lemma}

$Proof.$ We will need the following identities:

(A) \quad $m_{-1}m_0=[j+1]m_0m_{j-1}-[j]m_0m_j$,

(B) \quad $m_jm_{k+1}=[k+1]m_jm_{1}-[k]m_jm_{0}$,

(C) \quad $m_1m_j-m_0m_{j-1}=m_{k+1}m_{j-k}-m_km_{j-1-k}$.
\vskip .2cm

Claim $(A)$ and $(B)$ are easily proved by induction on $j$ (for $(A)$)
and on $k$ (for $(B)$), using Eq \ref{minduction}. Claim $(C)$ is shown
by induction on $k$, applying the two formulas in Eq \ref{minduction} to $m_{k+1}$ and to $m_{j-1-k}$.
Applying $(A)$ and $(B)$ for $j=k$, we obtain
$$
m_{-1}m_0-m_jm_{j+1}\ = [j+1](m_0m_{j-1}-m_jm_1).$$
Hence if $j=2k$, we obtain from $(C)$ that
$$m_0m_{j-1}-m_jm_1= m_km_{k-1}-m_{k+1}m_k=m_k(m_{k+1}-m_{k-1}).$$
Hence $\al_{j,j+1}=0$ implies that either $m_k=0$ or $m_{k+1}-m_{k-1}=0$, which implies $m_k=[0]_+$. Similarly, if $j=2k+1$, we have
$$m_0m_{j-1}-m_jm_1= m_k^2-m_{k+1}^2=(m_{k+1}+m_k)(m_{k+1}-m_k).$$
Observing that $m_{k+1}=\frac{1}{2}([2]m_k\pm\sqrt{(q-q^{-1})^2m_k^2+4}$, we deduce
from $m_{k+1}=\pm m_k$ that
$$(2\pm [2])m_k=\pm \sqrt{(q-q^{-1})^2m_k^2+4}.$$
This can be transformed to
$$m_k^2=\frac{1}{2\pm (q+q^{-1})}=\left[\frac{q^{1/2}\pm q^{-1/2}}{q-q^{-1}}\right]^2,$$
from which one deduces $m_k=\pm [1/2]_+$ or $m_k=[\pm 1/2]$, depending on
the choice of sign and the choice of square root of $q$. $\Box$

\subsection{Finite-dimensional modules for generic $q$}

\begin{theorem}\label{Vermafin}
We assume $q$ not to be a root of unity except possibly for $q=\pm 1$.
The Verma module $V_m$ of $U'_q\so_3$ has a finite-dimensional quotient
if

(a) $m=[\la]$, $m_1=[\la-1]$ or $m=[-\la]$, $m_1=[1-\la]$, for $\la\in\frac{1}{2}\Z$, $\la>0$, or

(b) $m=\pm [\la]_+$, $m_1=\pm[\la-1]_+$ for $\la\in\frac{1}{2}+\Z$, $\la>0$, with matching signs.

In both cases, the largest finite-dimensional
quotient has dimension $2\la+1$. It is simple
in case (a), and the direct sum $L_+\oplus L_-$ of two simple submodules of dimension $\la+1/2$ each in case (b). These submodules correspond
to the eigenspaces of $B_2$ with eigenvalues $[\la-j]_+$  respectively with eigenvalues
$-[\la-j]_+$, $j=0,1,\ ...\ \la-1/2$.
\end{theorem}

$Proof.$ By Proposition \ref{weightbasis}, the module $V_m$ has a
basis of weight vectors $v_j$ in the cases listed here.
It follows from a standard argument that then
also any submodule of $V_m$ would have a basis of weight
vectors. If $v_j$ is the vector of highest weight of such a submodule,
it follows that
$\al_{j-1,j}= 0$. One checks from Corollary \ref{ali}
that this happens for $j=2\la+1$ in both cases (a) and (b),
and that
$\al_{j-1,j}\neq 0$ for all $j\neq 2\la +1$. This shows that we have
finite dimensional quotients of dimension $2\la+1$, and any other
finite-dimensional quotient would have to be smaller.
If $m=[\la]$, all eigenvalues of $B_1$ are distinct in that quotient $L$.
As $\al_{j,j+1}\neq 0$ for $0\leq j<2\la$, one checks easily
that $L$ is simple.

If $m=[\la]_+$, the vectors $v_j$ and $v_{2\la-j}$ have the same eigenvalue
for $B_1$, for $0\leq j\leq 2\la$. We leave it to the reader to check
that after a suitable rescaling of these basis vectors, $B_2$ is represented
by a symmetric matrix $A$ such that $a_{j-1,j}=a_{j,j-1}=\sqrt{\al_{j-1,j}}$,
and $a_{i,j}=0$ if $|i-j|\neq 1$.
Hence $V_m/M$ is equal to the direct sum $L_+\oplus L_-$, where the modules
$L_\pm$ have bases $v_j\pm v_{2\la-j}$, $0\leq j\leq \la-\frac{1}{2}$.
Now again the eigenvalues of $B_1$ are mutually distinct on each of these
submodules, and one shows irreducibility as for case (a).
We will prove the last statement in the proof of Theorem \ref{irrepnon}. $\Box$

\subsection{The case $m=[\la]_+$, with $\la\in\Z$} As noted in
Propostion \ref{weightbasis},
we do not have a basis of weight vectors for the Verma module in this case.
In the following let $r$ be a positive integer.
We then define, for generic highest weight $m$, the vectors $v'_{r+j}$, $1\leq j\leq r$ by
$$v'_{r+j}=v_{r+j}+[\prod_{i=1}^j\al_{r-i,r-i+1}] v_{r-j},\quad 1\leq j\leq r.$$

\begin{lemma}\label{welldefined} The vectors $v'_{r+j}$, $1\leq j\leq r$
 are well-defined also if $m=[r]_+$.
\end{lemma}

$Proof.$ We prove the claim by induction on $j$. Observe that $v'_{r+1}=B_2v_r$, by Proposition \ref{weightbasis}.
Using Eq \ref{prop32}, one obtains
$$B_2v'_{r+1}=B_2^2v_r=v'_{r+2}+(\al_{r,r+1}+\al_{r-1,r})v_r,$$
where one checks that the coefficient of $v_r$ is well-defined
also for $m=[r]_+$.
Using the definitions, we similarly obtain,
for $j> 1$,
$$B_2v'_{r+j}= v'_{r+j+1} + \al_{r+j-1,r+j}v'_{r+j-1}+ (\al_{r-j,r-j+1}-\al_{r+j-1,r+j})[\prod_{i=1}^{j-1}\al_{r-i-1,r-i}]\ v_{r-j+1}.$$
It can then be shown by a direct calculation that the scalar of $v_{r-j+1}$
is also well-defined for $m=[r]_+$; the pole at $\al_{r-1,r}$
cancels with the zero at $\al_{r-j,r-j+1}-\al_{r+j-1,r+j}$ for $\la=r$.
Hence we can express $v'_{j+r-1}$ by an expression of vectors which are also
well-defined at $\la=r$, using the last formula. $\Box$

\begin{lemma}\label{mplus} The Verma module $V_m$ has no nontrivial submodule
for $m=[r]_+$, with $r$ a nonnegative integer.
\end{lemma}

$Proof.$ \ignore{ If $m=\pm [0]_+$, one calculates from Corollary \ref{ali}, (b)
that $\al_{j-1,j}=-2/(q-q^{-1})^2\neq 0$
for all $j>0$.}  By the previous lemma,
we have a basis for $V_m$ consisting of the vectors $v_j$, with $0\leq j\leq r$
or $j\geq 2r+1$, and the vectors $v'_{r+j}$, $1\leq j\leq r$. Now observe that
\begin{align}
B_2v_{2r+1}\ &=\ v_{2r+2}+\al_{2r,2r+1}v_{2r}\cr
&=\ v_{2r+2}+\al_{2r,2r+1}v'_{2r} - \al_{2r,2r+1}[\prod_{i=1}^{r-1}\al_{r-i-1,r-i}]\ v_0\cr
&=\ v_{2r+2}-\beta v_0 \quad {\rm for\ } m=[r]_+,\cr
\end{align}
where $\al_{2r,2r+1}=0$ and where one checks as in the proof of Lemma \ref{welldefined}
that for $m=[r]_+$ the scalar $\beta$ is well-defined and non-zero.
Assume now $M$ is a submodule of $V_m$, and $v\in M$ a nonzero vector.
Writing it as a sum of generalized eigenvectors of $B_1$, we can also
assume that each of the generalized eigenvectors is in $M$. Applying
$B_1$ to a generalized eigenvector, if necessary, we can also assume that
the corresponding eigenvector is contained in $M$, i.e. a vector $v_j$
with $j\leq r$ or $j\geq 2r$. In both cases, we can show that then also
$v_{2r+1}$ is in $M$, and hence also the generating vector $v_0$,
by the calculation above. $\Box$

\subsection{Irreducible $U'_q\so_3$-modules}\label{irrepmod}
 We can now classify
irreducible $U'_q\so_3$-modules for $q$ not a root of unity.
These results have been obtained
before, see \cite{IK} and references there, or \cite{HP}.

\begin{theorem}\label{irrepnon} Let $q$ not be a root of unity, except
possibly for $q=\pm 1$. Then the following gives a complete list of
finite dimensional simple modules $L$ of $U'_q\so_3$, up to isomorphism:

(a) The module $L$ has highest weight $m=[\la]$ with $m_1=[\la-1]$ and $\la\in\frac{1}{2}\Z$,
$\la\geq 0$. In this case, $L$ has dimension $2\la+1$, and it is uniquely
determined by its highest weight.

(b) The module $L$ has highest weight $m=\pm [\la]_+$ with $m_1=\pm [\la-1]_+$ (matching signs) and $\la\in\frac{1}{2}+\Z$,
$\la\geq 0$. In each of these cases,
$L$ has dimension $\la+\frac{1}{2}$, and there
are exactly two non-equivalent modules with these properties. They are
related by the outer automorphism $B_2\mapsto -B_2$, $B_1\mapsto B_1$.

(c) For any positive integer $k$, there are exactly 5 non-equivalent simple $U'_q\so_3$ modules of dimension $k$,
one in case (a) and four in case (b).
Those four cases are related by possible sign changes $B_j\mapsto -B_j$,
$j=1,2$.

(d) In every simple finite dimensional representation $L$ of $U_q'\so_3$, $B_1$ is conjugate
to $ B_2$ or to $-B_2$.
\end{theorem}

$Proof.$ We have shown in Lemma \ref{mplus} that the Verma module $V_m$
does not have a finite-dimensional quotient if $m=[\la]_+$ for
$\la$ a positive integer. In all other cases, the Verma module $V_m$
has a basis of weight vectors, by Proposition \ref{weightbasis}.
We can have a nontrivial submodule in $V_m$ only if $\al_{j-1,j}=0$
for some $j$, see e.g. the proof of Theorem \ref{Vermafin}.
It follows from Lemma \ref{finitelem} that this is possible only
if the weights are of the form $[\mu]$ with all $\mu\in \Z$ or all
$\mu\in \frac{1}{2}+\Z$, or all weights are of the form $\pm [\mu]_+$,
with $\mu\in \frac{1}{2}+\Z$ (see also Lemma \ref{mplus}).
Excluding the cases with $v_0$ being a lowest-weight vector
in the usual sense, we are left with the cases already listed in
Theorem \ref{Vermafin}.

Now observe that for $m=[\la]_+$, $B_1$ has trace
$\sum_{j=0}^{\la-1/2}[j]_+\neq 0$ on $L_\pm$ (see Theorem \ref{Vermafin}), 
and $B_2$ has only one nonzero
diagonal element, which is equal to $\pm a_{r-1,r}$, where $\la= r+\frac{1}{2}$.
This shows that $L_\pm$ are two non-equivalent
$U'_q\so_3$ modules. On the other hand, if $m=[\la]$, both $B_1$ and $B_2$ are
represented by matrices with trace equal to 0. This proves (c),
 as we can also replace $B_1$ by $-B_1$ in case (b);
changing the sign of the generators in (a) does not change the weights,
and hence must produce an equivalent module. Part(d) now follows
from the fact that constructing modules using the weight spaces of $B_2$
also only gives us 5 non-equivalent modules. By the trace argument above,
$B_1$ being in case (a) is equivalent to $B_2$ being in case (a).
Hence, as there is only one irreducible module in that case, $B_1$ and $B_2$
must have the same eigenvalues. The same argument also works for case (b),
up to possible sign changes. Finally, statements (d) and (c)
and the arguments in this paragraph also
imply the last statement of Theorem \ref{Vermafin}. $\Box$

\subsection{Examples} We now want to apply our results so far for
the representations of $U'_q\so_3$ on the $\ell$-dimensional module $V$
defined in Eq \ref{qtorus1} and \ref{qtorus2} for $q$ a primitive
$\ell$-th root of unity. For $\ell$ even, this has already more or less
appeared in \cite{RW}, Lemma 4.3.

\begin{lemma}\label{qtorusso3} (a) If $\ell$ is even, the module $V$ is
a direct sum of two irreducible modules of dimensions $\ell/2\pm 1$
with highest weights $[\ell/2]$ and $[\ell/2-1]$ for
Eq \ref{qtorus1}.

(b) If $\ell$ is odd, the module $V$ is the direct sum of two modules
of $U'_{q}\so_3$ of dimensions $(\ell\pm 1)/2$ and highest weights $-[\ell/2]_+$ and $-[\ell/2-1]_+$ for Eq \ref{qtorus1}.

(c) We only have highest weight modules for representations
given by Eq \ref{qtorus2} if $\ell$ is divisible by 4. In this case
they are equivalent to representations in (a).
\end{lemma}

$Proof.$
Let us do case (b) in some detail. Let $q^{1/2}$ be the square root
of $q$ for which $q^{\ell/2}=-1$. Using this and Eq \ref{qtorus1}, the eigenvalue
of $B_1$ for the $j$-th basis vector is
$$i\frac{q^j+q^{-j}}{q-q^{-1}}\ =\ -i\frac{q^{\ell/2-j}+q^{j-\ell/2}}{q-q^{-1}}.$$
This shows that the character of $V$
is the sum of the characters of the two highest weight modules as claimed.
Indeed, $v_0$ and $v_1-v_{\ell-1}$ are highest weight vectors for the given
weights.
It follows from Corollary \ref{ali} that the coefficients $\al_{j-1,j}$
in Proposition \ref{weightbasis} are well-defined and non-zero for the highest weights $[\ell/2]_+$ and $[\ell/2-1]_+$ also for our choice of $q$.
Hence these modules are irreducible. Part (a) is shown similarly and was already done in \cite{RW}.
Part (c) follows from the observation that we can not find a weight vector
$v\in V$ such that also $B_2v$ is a weight vector unless we can find
a weight $\pm 2i/(q-q^{-1})$. $\Box$

\section{Basic results for Verma modules of $U'_q\so_n$}

\subsection{Weights}  The following
notation will be convenient: If $\mu = (\mu_i)_i\in\Q^k$, we use
the notations $[\mu]$ for the vector $([\mu_i])$, and
$[\mu]_+$ for the vector $([\mu_i]_+)$, where, as usual
$[k]=(q^k-q^{-k})/(q-q^{-1})$ and $[k]_+=\sqrt{-1} (q^k+q^{-k})/(q-q^{-1})$.
If $k=r/s$ is rational, we assume that a choice for an $s$-th root of
$q$ has been made. The following lemma clarifies how generators which
are not in the Cartan algebra act on weight vectors. It will be convenient
to use the notation $\al_i^+=\ep_i+\ep_{i+1}$; here $\ep_i$ is as in
Section \ref{definitions}. Observe that $\al_i^+$ is a positive root.
Finally we also have the following

$Convention$: If $v$ is a weight vector of a Verma module which is not the
highest weight vector, we assume $B_{2i-1}v\neq \pm [0]_+ v$ for $1\leq i\leq \lfloor n/2\rfloor$.
 As we are
primarily interested in finite-dimensional modules, this is no serious restriction
in view of Lemma \ref{mplus}

\begin{lemma}\label{weightshift}
Let $v$ be a vector in a $U_q'\so_n$ module with weight $[\mu]$. Then

(a)\ $(B_{2i-1}-[\mu_i+1])(B_{2i-1}-[\mu_i-1])B_{2i}v=0$.

(b)\ The vector $B_{2i}v$ can be written as a linear combination of
at most four weight vectors, with weights $[\mu\pm \al_i]$ and
$[\mu\pm\al_i^+]$.

(c) The analogous statements also hold if every weight $[\nu]$ in (a) and (b)
is replaced by $[\nu]_+$.
\end{lemma}

$Proof.$ These are straightforward calculations. E.g. for (a) we have
\begin{align}
B_{2i-1}^2B_{2i}v\ &=\ ([2]B_{2i-1}B_{2i}B_{2i-1}-B_{2i}B_{2i-1}^2+B_{2i})v\cr
&=\ ([2][\mu_i]B_{2i-1}-([\mu_i]^2-1))B_{2i}v.
\end{align}
We now get the claimed factorization in (a) using the identities
$[2][\mu_i]=[\mu_i+1]+[\mu_i-1]$ and $[\mu_i]^2-1=[\mu_i+1][\mu_i-1]$.
For part (b)  observe that a similar calculation also holds
with $2i-1$ replaced by $2i+1$. As $B_{2i+1}$ commutes with $B_{2i-1}$, we obtain
weight vectors of the form
$(B_{2i+1}-[\mu_{i+1}\pm 1])(B_{2i-1}-[\mu_{i}\pm 1])B_{2i}v$, for any choice of signs for
which the expression is not equal to zero.  Part (c)
can be proved in the same way as parts (a) and (b). $\Box$

\subsection{Applications to general $U'_q\so_n$ modules}

\begin{corollary}\label{finiterestriction} Let $V$ be a finite-dimensional simple $U_q'\so_n$
module, for $q$ not a root of  unity except for $q=\pm 1$.
Then we have, using the notations of Section \ref{definitions} that

(a) all its weights are of the form $[\mu]$ with all $\mu_i\in\Z$ or
all $\mu_i\in\frac{1}{2}+\Z$, or

(b) all its weights are of the form $(\pm[\mu_i]_+)$ with all $\mu_i\in\frac{1}{2}+\Z$ (where $q\neq \pm 1$).
\end{corollary}

$Proof.$ This follows for $\mu_1$ from our results for $U'_q\so_3$.
Let $U'_q\so_3(i)$ be the subalgebra of $U'_q\so_n$ generated by $B_i$ and $B_{i+1}$;
it is isomorphic to $U'_q\so_3$.
Then we can show by induction on $i$, using Theorem \ref{irrepnon},(d)
that also the eigenvalues of $B_{i+1}$ have to be the same as the ones
of $B_i$, up to a possible sign change. $\Box$

\begin{definition}\label{domweights} We say that a weight $m$ of
$U'_q\so_n$ is regularly dominant if $m=[\la]$ or $m=[\la]_+$
with $\la$ a dominant integral weight of $\so_n$ with all $\la_i\in\frac{1}{2}
+\Z$ if $m=[\la]_+$
\end{definition}

\begin{proposition}\label{highestweights} Let $q$ not be a root of unity.
Then any simple finite dimensional $U'_q\so_n$-module $V$ is a quotient
of a Verma module with highest weight $m$ being regularly dominant,
after possible sign changes $B_{2i-1}\mapsto \pm B_{2i-1}$,
for $1\leq 2i-1\leq n$.
\end{proposition}

$Proof.$ Corollary \ref{finiterestriction} already associates the weights
of $V$ to weights in the classical case.
If we are in case (b),
we can assume that our module has at least one
weight of the form $[\mu]_+$ with all $\mu_i>0$,
after possible sign changes $B_{2i-1}\mapsto
\pm B_{2i-1}$.  As our module is finite dimensional, it must have a
highest weight. The rest of the proof now follows from usual standard arguments. If we are in case (a) of Corollary \ref{finiterestriction}, the same arguments
apply, without having to worry about possible sign changes. $\Box$

\begin{definition}\label{standardhw} We call $V_m=U'_q\so_n/I_m$
a standard Verma module
if $I_m=I_{m,n}$ is as before Proposition \ref{Vermaprop} with $m$ being
a regularly dominant weight, and with $n_i=[\la_i-1]$ respectively
$n_i=[\la_i-1]_+$ if $m_i=[\la_i]$ respectively $m_i=[\la_i]_+$.
\end{definition}

\begin{remark}\label{standardremark} If $q=1$ and $m=[\la]=\la$,
a standard Verma module will be the usual Verma module
with highest weight $\la$. As we do not have
explicit raising and lowering operators in the algebra $U'_q\so_n$
any vector whose weight only differs by sign changes in some of its
coordinates from the highest weight could also be made into a highest
weight vector of $U'_q\so_n$ according to the definition before Proposition \ref{Vermaprop}. E.g. the lowest weight vector of a simple finite
dimensional $U'_q\so_3$-module would also be a highest weight vector.
We can also carry out the theory of Verma modules for such weights.
As this would only lead to more tedious notations, without any new
results, we will restrict ourselves to standard Verma modules in the following.
\end{remark}

\subsection{Classical case $q=1$} In order to show linear independence of the spanning sets in Proposition \ref{Vermaprop}
we will appeal to the classical case at $q=1$. Using matrix units $E_{ij}$, we can identify $B_j$ with $\sqrt{-1}(E_{j,j+1}-E_{j+1,j})$.
It can be easily checked by explicit matrix calculations that we get
root vectors $X_\al$ of $\so_n$ as follows, where, for simplicity, we write $ad_i$ for $ad_{B_i}$, the adjoint operation by $B_i$:

1. If $\al=\pm \ep_i\pm \ep_j$, we define
$X_\al=(1\pm ad_{2i-1})(1\pm ad_{2j-1})ad_{2i}ad_{2i+1}\ ...\ ad_{2j-3}B_{2j-2}.$

2. If $\al=\pm \ep_i$ is a short root for $\so_n$, with $n$ odd, we have the root vector
$X_\al=(1\pm ad_{2i-1})ad_{2i}ad_{2i+1}\ ...\ ad_{n-2}(B_{n-1}).$

\ni
It will be convenient to use the following elements for describing a basis for the Verma module:
$$Y_{-\ep_j}=X_{\ep_j}+X_{-\ep_j}=2\ ad_{2j}ad_{2j+1}\ ...\ ad_{n-2}(B_{n-1}),$$
$$Y_{-\ep_i+\ep_j}=X_{\ep_i+\ep_j}+X_{\ep_i-\ep_j}+X_{-\ep_i+\ep_j}+X_{-\ep_i-\ep_j}
=4\ ad_{2i-1}ad_{2i+1}\ ...\ ad_{2j-3}(B_{2j-2}),$$

$$Y_{-\ep_i-\ep_j}=X_{\ep_i+\ep_j}-X_{\ep_i-\ep_j}-X_{-\ep_i+\ep_j}+X_{-\ep_i-\ep_j}
=4\ ad_{2i}ad_{2i+1}\ ...\ ad_{2j-3}(B_{2j-2}).$$
We will use the following well-known simple observation in the proof
of the following lemma:  Assume $(Y_i)_i$ is a set of elements in a semisimple Lie algebra $\g$
whose residues mod $\b_+$ form a basis of $\g/\b_+$; here $\b_+$ is the Borel algebra spanned by the Cartan
algebra and the root vectors corresponding to positive roots. Then
the Verma module $\Vl$ has a basis
\begin{equation}\label{Vermabasisequ}
\prod_i Y_i^{k_i}v_\la,\quad k_i\geq 0,
\end{equation}
where $v_\la$ is the highest weight vector. Indeed, this is well-known if
the $Y_i=X_i\in \n_-$ form a basis of the nilpotent part spanned by weight vectors corresponding to negative roots.
In the general case, we can write
$Y_i=X_i+Z_i$, with $X_i\in\n_-$ and $Z_i\in\b_+$. It is well-known that
the Verma module $\Vl$ has a filtration whose $r$-th subspace $\Vl(r)$ is spanned by all products of $\leq r$ generators
$X_i\in\n_-$ modulo $I_m$, and that $Z\Vl(r)\subset \Vl(r)$ for all $Z\in\b_+$.
It now follows easily that we get the same filtration also in terms of the generators $Y_i$.

\begin{lemma}\label{classic} Let $V_m$ be a standard Verma module, see Definition
\ref{standardhw}, with $m=[\la]$ for an integral dominant weight $\la$.
Then the spanning set in Proposition \ref{Vermaprop}, (b) is linearly independent for $q=1$.
\end{lemma}

$Proof.$ 
We  use the observation before this lemma for the ordering as defined before
Lemma \ref{wordord}. First observe that with respect to this
ordering we have, for $i<j$, 
$$ad_{i}ad_{i+1}\ ...\ ad_{j-1}(B_{j})\ =
\ (-1)^{j-i}B_jB_{j-1}\ ...\ B_{i+1}B_i\ +\ lower\ terms.$$
Hence, mapping $Y_\al$, with $\al$ a negative root, to its highest term
in its expansion of products of $B_m$'s, we obtain a bijection
$$\Phi: \{ Y_\al, \al<0\}\quad \leftrightarrow\quad  \{ B_{j,i}, 1\leq i<j\leq n, i\ {\rm even}\}.$$
As the ordering is preserved by multiplication, we obtain in general that
$$\Phi(\prod Y_\al^{k(\al)})= \prod \Phi(Y_\al)^{k(\al)} + lower\ terms.$$
Hence we get a triangular transformation matrix between the spanning
set $\B$ of Proposition \ref{Vermaprop},(b) and the basis in Eq \ref{Vermabasisequ}. $\Box$

\subsection{Linear independence} We will prove linear independence of the set $\B$  in 
Proposition \ref{Vermaprop}, (b) for standard Verma modules
by appealing to the already developed representation
theory of $U_q'\so_n$, see Section \ref{knownrep}. We give
the standard arguments in the proof of the following theorem for the
reader's convenience. We will use the following well-known
result, which follows
from a more precise result by
Harish-Chandra (see e.g. \cite{Va}, Theorem 4.7.3).

\begin{proposition}\label{injectweight} Let $q=1$ and
let $\la$ be an integral dominant
weight, with $V_\la$ and $L_\la$ the Verma module and the irreducible module
with highest weight $\la$ respectively. Let $\mu=\la-\sum_i r_i\al_i$
for some non-negative integers $r_i$, and let $\r=(r_i)_i$. We denote
by $V_\la[\r]$ and $L_\la[\r]$ the span of weight vectors in the respective
modules with weights of the form $\la-\sum \tilde r_i\al_i$, with $0\leq \tilde r_i\leq r_i$. Then, for given $\r$, the restriction of the
canonical map $\Phi: V_\la\twoheadrightarrow L_\la$ to $V_\la[\r]$ is injective for all but finitely many dominant integral weights $\la$.
\end{proposition}

$Proof.$ We give the short proof for the reader's convenience.
By Harish-Chandra's theorem (see e.g. \cite{Va}, Theorem 4.7.3),
the kernel of $\Phi$ is spanned by submodules of $V_\la$ with highest
weights $s_i.\la=\la-((\la,\al_i)+1)\al_i$. As all the weights $\mu$
of such a submodule satisfy $\mu\leq s_i.\la$,
it follows that no weight vector
of $V_\la$ with weight $\la-\sum \tilde r_i\al_i$ for which
$\tilde r_i\leq (\la,\al_i)$
for all $i$
can be in the kernel of $\Phi$. Hence, if $r_i=(\la,\al_i)$,
none of the weights of $V_\la[\r]$ can be in the kernel of $\Phi$. $\Box$

\begin{theorem}\label{linindep}
The set $\B$
in Proposition \ref{Vermaprop} (b) is a basis for the Verma module
$V_{m,\tilde n}$, provided that $m_i^2-[2]m_i\tilde n_i+\tilde n_i^2=1$ for $1\leq i\leq \lfloor n/2\rfloor$. 
This includes, in particular, standard Verma modules, see Def \ref{standardhw}.
\end{theorem}

$Proof.$ Let $\B$ be as in Proposition \ref{Vermaprop}, (b),
and let $S=\{ M_j$, $1\leq j\leq k\}$ be a finite subset of $\B$.
By Lemma \ref{classic}, their image in the Verma module $V_\la$
is linearly independent for any dominant weight $\la$ for $q=1$.
By Lemma \ref{weightshift}, (b), there exists a vector $\r$ with non-negative
integer entries such that the image of $S$ is contained in $V_\la[\r]$.
By Proposition \ref{injectweight}, linear independence even holds for the
image of $S$ in $L_\la$, for all but finitely many dominant integral weights. 
As discussed in Section \ref{defrep}, the representations $L_\la$
are also well-defined for generic $q$. In particular, the image
of $S$ is also linearly independent in $L_\la$ for all but finitely
many dominant integral weights, and for generic $q$.

We now want to prove linear independence for general $q$
for Verma modules $V_{m,\tilde n}$ with $m_i=[\la_i]$ and $\tilde n_i=[\la_i-1]$;
for simplicity of notation, we will write $V_\la$ for $V_{m,\tilde n}$.
We define a linear action
of the generators $B_i$ on a vector space whose basis elements
are labeled by the monomials in $\B$, using the relations in
the proof of Proposition \ref{Vermaprop} and its preceding lemmas.
Fix a monomial $M\in\B$
and generators $B_i$ and $B_{i+1}$. Apply each monomial which appears
in relation \ref{basicrel} to $M$ and expand this as a linear combination
of elements in $\B$. Let $S$ be the finite subset of $\B$ consisting of
all those elements appearing in these expansions.
It follows from the discussion in the previous paragraph that
the image of $S$ is linearly independent in $L_\la$ for all but finitely
many dominant integral weights $\la$. As $L_\la$ is a $U'_q\so_n$-module,
it follows that
the relation \ref{basicrel} holds, if applied to $M$, for all but
finitely many $\la$'s. Observe that
the matrix coefficients for $B_i$ and $B_{i+1}$
are Laurent polynomials in the variables $q$ and $q^{\la_i}$,
and the elements  $(q-q^{-1})^{-1}$ and $(q+q^{-1})^{-1}$ (see Proposition \ref{Vermaprop} and
the lemmas used for it). Hence the relation holds for a Zariski-dense
set of values $m_i=[\la_i]$
and $\tilde n_i=[\la_i- 1]$, for highest weight $m=[\la]$; this implies it holds
for such $m$ and $\tilde n$ for any choice of values for $\la$.
It follows that we have constructed a highest weight representation
of $U'_q\so_n$ with a basis labeled by the elements of $\B$.
Hence $\B$ itself must
be linearly independent for these choices of parameters.

To conclude the proof, we also need to show the claim for all possible
choices of $\tilde n_i$, for given $m$. We have already seen that if $m_i=[\la_i]$,
the two possible solutions for $\tilde n_i$ are $[\la_i\pm 1]$ (see the examples
before Lemma \ref{mbasic}). We can mimic the proof above by using
integral weights of the form $(\la_i)$, where we choose negative entries
for those coordinates $i$ for which we want $\tilde n_i=[\la_i+1]$. By defining
the Weyl chamber of this $\lambda$ to be the dominant Weyl chamber,
we again obtain a finite dimensional quotient; it is isomorphic to
$L_{|\la|}$, where $|\la|$ is the weight in the usual dominant Weyl
chamber which is conjugate to $\la$. We can now show as in the previous
paragraph that $\B$ is linearly independent for any $\la$, with
$m_i=[\la_i]$ and $\tilde n_i=[\la_i\pm 1]$ for our given choice of signs.
This concludes the proof of linear independence of $\B$ for any $m$ and
$\tilde n$ as given in the statement. $\Box$

\begin{corollary}\label{sameweights}
If $m=[\la]$ and $V_m$ is a standard Verma module 
for $q$ not a root of unity (see Definition \ref{standardhw}), 
we have the `same' weight multiplicities for $V_m$
as in the classical case $q=1$, i.e.
the weight $[\mu]$ has the same multiplicity
in $V_m$ as the weight $\mu$ has in $V_\la$ for $q=1$.
\end{corollary}

$Proof.$ If $V_{m}$ is a standard
Verma module with $m=[\la]$, then its weights
are of the form $[\la-\om]$, with $\om$ in the root lattice, by Lemma \ref{weightshift}. This is also true for its
finite-dimensional subspaces
$W=W(0)=V_{m}[[\s ]]$, as defined before Corollary \ref{Vermapropcor}.
The multiplicities of the (generalized) weight spaces are obtained
from the characteristic polynomials of the $B_i$'s acting on
$W$ as follows: Fix a weight $\mu$. Let $W(1)$ be the eigenspace of $B_1$
for $[\mu_1]$. Assuming we have constructed $W(i)$, define $W(i+1)$ to be the eigenspace of
$B_{2i+1}$ in $W(i)$ for eigenvalue $[\mu_{i+1}]$, if $2i+1<n$.
Then the weight space $W[\mu]$ is equal to $W(\lfloor (n-1)/2\rfloor)$.
The multiplicity of the eigenvalue $[\mu_i]$ in the characteristic polynomial of $B_{2i-1}$ 
on $W(i-1)$ does not change with $q$ as long as it is not a root of unity $\neq 1$.
Hence the dimensions of the $W(i)$ stay constant.
It follows that the statement of the corollary
is true for these subspaces.
For any given weight $\mu$, we can find a
vector $\s$ with sufficiently large coordinates such that
the multiplicity of $\mu$ in $V_{m}[[\s ]]$ coincides
with the one of $V_{m}$. The claim follows from this. $\Box$

\begin{remark}\label{weightmult}
If $m=[\la]_+$ a slight subtlety occurs stemming from the
fact that $[r]_+=[-r]_+$ for any rational number $r$.
However, if the $\la_i$s are not
integers, we can inductively decide for each eigenvector of $B_{2i-1}$
whether it has formal eigenvalue $[\nu_i]_+$ or $[-\nu_i]_+$, by Lemma \ref{weightshift}, (a). 
Indeed, if $v$ is a weight vector with weight $[\mu]_+$, and if we write $B_{2i}v$ as a linear combination
of weight vectors,
we would only have difficulties in determining their formal weights if $\mu_i=0$
or $\mu_{i+1}=0$. Hence we can and will define distinct weight spaces
$V_m[\mu]$ and $V_m[\tilde\mu]$ for weights $\mu\neq\tilde\mu$ for which
$|\mu_i|=|\tilde\mu_i|$ for all $i$, even though the corresponding weights $[\mu]_+$
and $[\tilde\mu]_+$ coincide. We can then again show similarly as in the proof of Coroallry \ref{sameweights}
that the multiplicty of the weight $[\mu]_+$ in $V_m$ for $m=[\lambda]_+$ is
the same as the multuplicity of the weight $\mu$ in $V_\la$ for $q=1$.
\end{remark}

\section{Representations of $U'_q\so_n,\ n\geq 4$}

\subsection{Preliminaries}\label{so4repsec}
We first need some explicit results for $U'_q\so_4$.
It follows from  Proposition \ref{highestweights}
that any finite-dimensional $U'_q\so_4$ module is a quotient of a standard
Verma module $V_m$.
Moreover, the weights of $V_m$
are of the form $[\mu]=[\la-r_1\al_1-r_2\al_2]$ for $m=[\la]$ or of the form
$[\mu]_+=[\la-r_1\al_1-r_2\al_2]_+$ for $m=[\la]_+$,
where $r_1,r_2\geq 0$ and $\al_1=\ep_1-\ep_2$ and $\al_2=\ep_1+\ep_2$.
Using the well-known isomorphism $\so_4\cong\sl_2\oplus\sl_2$, we see
each weight in the Verma module has multiplicity 1 in the classical case.
So we have a basis
$(v_\mu)$, with $\mu$ running through the set as described above for any
standard Verma module,
uniquely determined up to scalar multiples. In view of Theorem \ref{linindep}, we can use the same notation also for a basis of weight
vectors for a standard Verma module $V_m$ of $U'_q\so_4$. In the following, it will be
convenient to use the notation and easily proved identity
\begin{equation}\label{curl}
\{k\}=q^k+q^{-k}=[k+1]-[k-1].
\end{equation}

\begin{lemma}\label{so4normalization}
Let $\mu=\la-r_1\al_1-r_2\al_2$. The following statements are given
for $m=[\la]$. They similarly hold for $m=[\la]_+$ after replacing
any numbers $[a]$ by $[a]_+$.

(a) There exists a unique weight vector $v_\mu$ of weight $\mu$ of the form
$$v_\mu=((B_3B_2)^{r_1+r_2}+ lower\ terms)v_\la.$$

(b) If $\mu=\la-r_i\al_i$, $v_\mu$ can be defined inductively by
$v_\mu=(B_3-[\la_2\pm (r_i-2)])B_2v_{\la-(r_i-1)\al_i}$,
where we have a plus sign for $i=1$ and a minus sign
for $i=2$ respectively. 
In particular, there exists
a polynomial $P_{\la,r,i}$ in two non-commuting variables such that
$v_{\la-r\al_i}=P_{\la,r,i}(B_2,B_3)v_\la$.
\end{lemma}

$Proof.$ We do the proof here for $m=[\la]$.
It goes by induction on $r_1+r_2$, which is trivially
true for $r_1+r_2=0$. 
It follows from Lemma \ref{weightshift} that $B_2v_\mu$ is
a linear combination of weight vectors of weights $(\mu_1\pm 1,\mu_2\pm 1)$.
Hence $(B_3-[\mu_2-1])B_2v_\mu$ is a linear combination of two
weight vectors of weights $(\mu_1\pm 1,\mu_2+1)$. Subtracting a
suitable multiple of $v_{(\mu_1-1,\mu_2+1)}$, we obtain a vector of
weight $(\mu+1,\mu_2+1)$. The statement about the leading
term follows by induction assumption.
Part (b) is a special case of part (a). Here it is easier to write
down an explicit formula due to the fact that
$\la-r\al_1+\al_2$ and $\la-r\al_2+\al_1$ are not weights for $V_m$. $\Box$

\begin{lemma}\label{so4rep} (see \cite{HKP}, Section VI)
Let $q$ not be a root of unity, let $m=[\la]$ and
let $(v_\mu)$ be the basis of weight vectors of $V_m$ as in
Lemma \ref{so4normalization}.
The action of the generators on
$v_\mu$,
with $\mu=\la-r_1\al_1-r_2\al_2$, is given by $B_1v_\mu=[\mu_1]v_\mu$, $B_3v_\mu=[\mu_2]v_\mu$ and
\begin{align}\label{so4repform}
B_2v_\mu=&
-\frac{\{\la_1-r_1+1\}\{\la_2+r_1\}[r_1][\la_1-\la_2-r_1+1]}
{\{\mu_1\}\{\mu_1+1\}\{\mu_2\}}\ v_{\mu+\al_1}\cr
& -\frac {\{\la_1-r_2+1\}\{\la_2-r_2\}[r_2][\la_1+\la_2-r_2+1]}
{\{\mu_1\}\{\mu_1+1\}\{\mu_2\}}\ v_{\mu+\al_2}\cr
& +\frac{1}{\{\mu_2\}}[\ v_{\mu-\al_1}-v_{\mu-\al_2}].
\end{align}
\end{lemma}

$Proof.$ Up to a renormalization of the basis vectors, this is essentially
just an extension of the results in e.g. \cite{HKP}, Section VI from
the finite dimensional modules to the
full Verma module. We give a brief outline here how this could be proved
similarly to our approach for $U'_q\so_3$ modules, see Proposition \ref{weightbasis}.
It follows from Lemma \ref{weightshift} that we can write $B_2v_\mu$
as
\begin{equation}\label{so4prelim}
B_2v_\mu= a_{\mu,\al_1}v_{\mu+\al_1}+a_{\mu,\al_2} v_{\mu+\al_2} + \frac{1}{\{\mu_2\}}v_{\mu-\al_1}-\frac{1}{\{\mu_2\}}v_{\mu-\al_2},
\end{equation}
for suitable scalars $a_{\mu,\al_i}$, $i=1,2$; here the scalars
for $v_{\mu-\al_i}$ are determined by the normalization of
our weight vectors as in Lemma \ref{so4normalization}.
The remaining scalars are calculated by induction on
$r_1+r_2$ by explicitly checking the relation
$B_2^2B_1-[2]B_2B_1B_2+B_1B_2^2=B_1$ at the vector $v_\mu$, with
$\mu$ as above.
Comparing the coefficients of the vectors $v_{\mu-\al_1+\al_2}$
and $v_{\mu+\al_1-\al_2}$, we obtain the recursion relations
$$\{\mu_1-1\}\{\mu_2+1\}a_{\mu-\al_1,\al_2}=\{\mu_1+1\}\{\mu_2\}a_{\mu,\al_2},$$

$$\{\mu_1-1\}\{\mu_2-1\}a_{\mu-\al_2,\al_1}=\{\mu_1+1\}\{\mu_2\}a_{\mu,\al_1}.$$
Moreover, comparing the coefficients of $v_\mu$ in the same relation,
one obtains
$$\frac{\{\mu_1-1\}}{\{\mu_2\}}(a_{\mu-\al_2,\al_2}-a_{\mu-\al_1,\al_1})
= -[\mu_1]+\{\mu_1+1\}(\frac{1}{\{\mu_2+1\}}a_{\mu,\al_2}-
\frac{1}{\{\mu_2-1\}}a_{\mu,\al_1}).$$
Finally, comparing the coefficients of $v_\mu$ after applying the relation
$B_2^2B_3-[2]B_2B_3B_2+B_3B_2^2=B_3$ to that vector, we also obtain
$$\frac{1}{\{\mu_2\}}(\{\mu_2-1\}a_{\mu-\al_2,\al_2}+
\{\mu_2+1\}a_{\mu-\al_1,\al_1})
= -[\mu_2]+a_{\mu,\al_1}+a_{\mu,\al_2}.$$
Using the last two recursion relations, one can calculate the coefficients
$a_{\la-r_1\al_1,\al_1}$ and $a_{\la-r_2\al_2,\al_2}$ by induction on
$r_1$ and $r_2$ respectively. Using the first two relations then
allows us to calculate the coefficients in general. $\Box$

\begin{corollary}\label{so4repp} If $q$ is not a root of
unity and $m=[\la]_+$, with $\la_i$ not a
positive integer for $i=1,2$, the Verma module $V_m$
has a basis of weight vectors and we obtain representations
as in the previous lemma by
$B_1v_\mu=[\mu_1]_+v_\mu$ and $B_3v_\mu=[\mu_2]_+v_\mu$ and by replacing
every occurrence of a factor
$\{m\}$ by $(q^m-q^{-m})$ in the definition of the action of $B_2$.
\end{corollary}

$Proof.$ The proof goes completely analogous to the one of Lemma \ref{so4rep}. $\Box$

\medskip

\subsection{Quotients of Verma modules for $U'_q\so_4$}
As usual, we define the dot action of the Weyl group $W$ on weights
by $w.\mu=w(\mu+\rho)-\rho$, where $\rho$ is half the sum of the positive roots. 
If $m=[\la]$ or $m=[\la]_+$,
we define $w.m$ to be equal to $[w.\la]$ and to $[w.\la]_+$ respectively.
If $\g=\so_4$ and
$s_1(\mu_1,\mu_2)=(\mu_2,\mu_1)$ and
$s_2(\mu_1,\mu_2)=(-\mu_2,-\mu_1)$, we obtain
$s_1.\la=\la-(r_1+1)\al_1$, $s_2.\la=\la-(r_2+1)\al_2$ and
$s_1s_2.\la=\la-(r_1+1)\al_1-(r_2+1)\al_2$, where $r_1=\la_1-\la_2$
and $r_2=\la_1+\la_2$.

\begin{proposition}\label{Vermaso4} Let $V_m$ be a standard Verma module 
(see Definition \ref{standardhw}) of
$U'_q\so_4$, with $m=[\la]$ or $m=[\la]_+$. We will write $w.\la$
for $w.m$ for simplicity of notation.

(a) The module $V_m$ has highest weight vectors 
$v_{s_i.\la}=P_{\la, (\la,\al_i)+1,i}(B_2,B_3)v_\la$
with weights
$s_i.m$, $i=1,2$, where $P_{\la, (\la,\al_i)+1,i}$ is as in
Lemma \ref{so4normalization}. Moreover,
each of these vectors generates a submodule of
$V_m$ which is isomorphic to the Verma module with the same weight.

(b) Let $M$ be a $U'_q\so_4$ module with highest weight $\la$.
Assume there exists a constant $K$ such that $|\mu_1-\mu_2|\leq K$
for each weight $\mu$ of $M$.  Then the weight spaces $M[\mu]$ and
$M[s_1(\mu)]$ have the same dimensions for each weight $\mu$ of $M$.
\end{proposition}

$Proof.$ Let $m=[\la]$. Observe that $s_i.\la=\la-((\la,\al_i)+1)\al_i$.
If $\mu=\la-r_i\al_i$, then Eq \ref{so4repform} simplifies to
\begin{equation}\label{so4simple}
B_2v_\mu= (-1)^{i-1}\frac{[r_i][(\la,\al_i)+1-r_i]}{\{\mu_1\} }v_{\mu+\al_i} +
\frac{1}{\{\mu_2\} }(v_{\mu-\al_1}-v_{\mu-\al_2}).
\end{equation}
It follows that the vector denoted by $v_{s_i.\la}$
is a highest weight vector. The proofs for $m=[\la]_+,
\la_i\in \frac{1}{2}+\Z$, $i=1,2$ go similarly.

By abuse
of notation, we will denote the submodules of $V_m$ generated by the
vector $v_{w.\la}$ by $V_{w.\la}$. By assumption, we have a
surjective map $\Phi$ from $V_m$ onto $M$. It follows again from
Lemma \ref{so4rep} that the only nontrivial submodule of
$V_{s_1.\la}$ is  $V_{s_1s_2.\la}$. We do know all the weights
of $V_{s_1.\la}$ and $V_{s_1.\la}/V_{s_1s_2.\la}$, from which
we see that these two modules do not satisfy the condition on the
weights as in the statement. Hence the kernel of $\Phi$ must contain
$V_{s_1.\la}$. We conclude that $M$ must be isomorphic to $V_\la/V_{s_1.\la}$
or to $V_\la/(V_{s_1.\la}+V_{s_2.\la})$, both of which satisfy the
claim. $\Box$

\subsection{Quotients of Verma modules for $U'_q\so_n$}

\begin{lemma}\label{generalideals}
Let $V_m$ be a standard Verma module for $U'_q\so_n$.
If $s_i$ is a simple reflection of the Weyl
group $W$ of $\so_n$, then $V_m$ has a submodule with highest
weight $s_i.m$.
\end{lemma}

$Proof.$ Let $U'_q\so_4(i)$ be the subalgebra of $U'_q\so_n$ which is
generated by $B_{2i}$ and $B_{2i\pm 1}$. By Proposition \ref{Vermaprop},
is isomorphic to $U'_q\so_4$.
We define $v_{s_i.\la}=P_{\la, (\la,\al_i)+1,1}(B_{2i},B_{2i+1})v_\la$, 
where $P_{\la, (\la,\al_i)+1,1}$
is the polynomial in Lemma \ref{so4normalization},
for $1\leq i\leq \frac{n}{2}-1$. If $i=n/2$ for $n$ even,
we define $v_{s_i.\la}=P_{\la, (\la,\al_{n/2})+1,2}(B_{n-2},B_{n-1})v_\la$, 
and if $n$ is odd
and $i=(n-1)/2$, we define $v_{s_i.\la}=P_{2\la_{(n-1)/2}+1}(B_{n-1})v_0$,
with the polynomial as in Corollary \ref{so3pol}.
We claim that each $v_{s_i.\la}$ is a highest weight vector.
If $j$ is odd, then the generator $B_j$ either commutes with $U'_q\so_4(i)$,
or it is an element of it. One deduces from this that the vector $v_{s_i.\la}$
is a weight vector. By the same argument, one easily checks the highest weight
property for all $B_{2j}$, except for $B_{2i+2}$. By Lemma \ref{weightshift},
$B_{2i+2}v_{s_i.\la}$
is a linear combination of vectors of weights $s_i.\la\pm \al_{i+1}$
or $s_i.\la\pm \al_{i+1}^+$. But as $s_i.\la=\la-r\al_i$ for a suitable
multiple $r$, the weights above can only be weights of the Verma module
if we have a minus sign at $\pm$. Hence $B_{2i+1}B_{2i+2}v_{s_i.\la}$
is a multiple of $B_{2i+2}v_{s_i.\la}$, i.e. the highest weight property
is satisfied. $\Box$

\begin{theorem}\label{character} Let $V_m$ be a standard Verma module,
(see Definition \ref{standardhw}) with $m=[\la]$ or $m=[\la]_+$.
 Let $I(\la)=\sum_i V_{s_i.\la}$,
where $V_{s_i.\la}$ is the highest weight module generated by $v_{s_i.\la}$,
as defined in Lemma \ref{generalideals}.
Then $V_m/I(\la)$ is a finite dimensional module whose character and,  in
particular, its dimension coincides with the one of the irreducible
$\so_n$ module with highest weight $\la$.
\end{theorem}

$Proof.$ We mimic the classical proof by showing that the Weyl group
permutes the weight spaces of the quotient module (see e.g \cite{Va}, Theorem 4.7.3
and preceding lemmas/theorems). Fix a simple reflection
$s_i$, where we assume $i\leq (n-2)/2$ for the moment.
Let $U$ be the subspace of $V_m/I(\la)$
consisting of vectors $u$ which generate a $U'_q\so_4(i)$ module $M(u)$
which is symmetric under $s_i$; by this we mean that the weight spaces
$M(u)[\mu]$ and $M(u)[s_i(\mu)]$ have the same dimensions, for all weights
$\mu$ of $M(u)$. It follows from Proposition \ref{Vermaso4} that $v_\la\in U$.

Next we claim that if $u\in U$, then so is $B_ju$, for $1\leq j<n$.
This is easy to see if $B_j$ commutes with $U'_q\so_4(i)$, and obvious
if $B_j\in U'_q\so_4(i)$. Hence we only need to worry about $B_{2i\pm 2}$.
Now observe that
$$U'_q\so_4(i)B_{2i+2}M(u)\subset M(u)+\sum_{j=2i-1}^{2i+2}B_jB_{j+1}\ ...\ B_{2i+2}M(u);$$
to see this, it suffices to show that the right hand side is a
$U'_q\so_4(i)$-module which contains $B_{2i+2}M(u)$. This is a consequence
of the relations proved in the lemmas before Proposition \ref{Vermaprop}.
As $M(u)$ is symmetric under $s_i$,
it follows that $|\mu_i-\mu_{i+1}|\leq K$
for a fixed constant, and any weight $\mu$ of $M(u)$.
But then it follows from the last inclusion that there also exists
a finite constant for the weights of $U'_q\so_4(i)B_{2i+2}M(u)$.
The proof for $B_{2i-2}$ goes similarly.
Hence this module is symmetric under $s_i$,
by Proposition \ref{Vermaso4}. It follows that $U=V_m/I(\la)$.
The proof for $i=n/2$
for $n$ even goes similarly, while the proof for $i=(n-1)/2$
is easier, only involving the corresponding argument in the setting of
$\so_3$.

It follows from the last paragraph that $V_m/I(\la)$ is symmetric with
respect to every simple reflection, and hence with respect to its Weyl
group. Hence the dimension of the weight space for $\mu$ is equal to
the dimension of the weight space belonging to its conjugate in the
dominant Weyl chamber. As our ideal does not contain any weights
in that region, by definition, this dimension coincides with the dimension
of the weight space in the Verma module. Hence the character of our
quotient module coincides with the character of the simple $\so_n$ module
with highest weight $\la$. $\Box$

\begin{corollary}\label{Weylmplus} If $m=[\la]_+$ with
$\la_i\in\frac{1}{2}+\Z$, then $V_m/I(\la)$ decomposes into the direct
sum of $2^{\lfloor (n-1)/2\rfloor}$ representations of equal dimension.
Any two of these can be made equivalent after suitable sign changes
$B_{2i}\to \pm B_{2i}$. Its weights are given by all weights $\mu$
of $V_m/I(\la)$ with $\mu_i>0$ for all $i$, with the same multiplicity
as in $V_m/I(\la)$, except possibly if $n$ is even, and $\la_{n/2}<0$.
In this case, the character is again the same as for $\bar\la$ which
coincides with $\la$ except for the last coordinate, where $\bar\la_{n/2}=-\la_{n/2}$.
\end{corollary}

$Proof.$ According to Theorem \ref{irrepnon}, (c) the simple $U'_q\so_3$
module
$L_m$, $m=[\la]_+$ with $\la\in \frac{1}{2}+\Z$ decomposes into a direct
sum $L_{m+}\oplus L_{m-}$, where $L_{m\pm}$ are the direct sum of
eigenspaces of $B_2$ corresponding to the eigenvalues $\pm [\la-j]_+$;
in the following, we also refer to these spaces $L_{m\pm}$ as eigenspaces
of $B_2$.
Using the same argument for the subalgebra of $U'_q\so_4$ generated by
$B_2$ and $B_3$, we similarly conclude that also $B_3$  leaves invariant
the eigenspaces $L_{m\pm}$ of $B_2$. We conclude that also $V_m/I([\la]_+)$
decomposes into the direct sum of two submodules for any admissible dominant
weight $\la$ of $U'_q\so_4$.

For the general case, we just construct analogous subspaces for
each $B_{2i}$. The statement about the character comes from the fact
that $[\mu_i]_+=[-\mu_i]_+$. $\Box$

\begin{remark}\label{irreducibilityrem} We have not shown here
that the quotients constructed in Theorem \ref{character}
for $m=[\la]$ respectively each of the $2^{\lfloor(n-1)/2\rfloor}$
components for $m=[\la]_+$ in Corollary \ref{Weylmplus} are irreducible.
This was shown in \cite{IK} in general, and for $m=[\la]$ also in
\cite{WSp}. It should be possible to prove this result directly in
the usual way. Casimir elements have been constructed in \cite{GI}
and \cite{NUW}, and the
scalar $c_\la$ via which it acts on an irreducible module with
highest weight $\la$ has been calculated in \cite{GI}. So it would suffice
to check that $c_\mu<c_\la$ whenever $\mu<\la$, i.e. $\la-\mu$ is a
sum of positive roots.
\end{remark}

\begin{theorem}\label{compclass} (see \cite{Kl}, \cite{IK})
Let $q$ not be a root of unity except for $q=\pm 1$, and let $V$ be
a simple finite-dimensional
$U'_q\so_n$-module. Then $V$ is either a classical simple
module with highest weight $[\la]$, where $\la$ is a dominant integral
weight for $\so_n$, or it is one of $2^{n-1}$ representations
belonging to the highest weight $[\la]_+$, where now $\la$
is a weight of $\so_n$ whose coordinates are
all in $\Z+\frac{1}{2}$ and positive.
In the second case, the dimension is $1/2^{\lfloor (n-1)/2\rfloor}$ times the dimension of the $\so_n$ module with highest
weight $\la$. The different modules can be obtained from each
other by multiplying
generators with even indices by $-1$.
\end{theorem}

$Proof.$ The restriction on the highest weights follows from
Proposition \ref{highestweights}. The existence of modules for such highest
weight modules follows from Theorem \ref{character} and Corollary
\ref{Weylmplus}. The dimensions of the finite-dimensional simple
modules for these highest weights was shown in \cite{IK} in general,
and for $m=[\la]$ also in \cite{WSp}. See also the discussion in
Remark \ref{irreducibilityrem}. $\Box$

\begin{remark}\label{Vermaresults}
The approach in \cite{IK} and its predecessors consisted of constructing
explicit matrix representations. Our explicit representations for $U'_q\so_3$
and $U'_q\so_4$ can be considered
as slight generalizations of special cases of their results. Unfortunately
their formulas become quite involved for larger $n$, and they are not well-defined for $q$ a root of unity.
\end{remark}

\section{Roots of Unity}

\subsection{Generic modules at roots of unity} In the following we assume
the representations to be defined over the complex numbers $\C$, with $q$ being
a primitive $\ell$-th root of unity.
Observe that in this case we have
\begin{equation}\label{rootshift}
[\la]_+\ =\ [\la+\ell/4].
\end{equation}
So, for $\ell$ even, we would not have to distinguish between classical
and non-classical representations. The following is an easy 
corollary of  Proposition \ref{Vermaprop} and  Theorem \ref{linindep}.

\begin{lemma}\label{linq} The set $\B$ in Proposition \ref{Vermaprop}
also defines a basis for a $U'_q\so_n$ module with highest weight 
as in Theorem \ref{linindep} for $q$ a root of unity. We
refer to this module as a Verma module.
\end{lemma}

$Proof.$ It follows from Proposition \ref{Vermaprop}
that the matrix coefficients of the action of a generator
$B_i$  with respect to the basis $\B$ are
in the ring $R[m,\tilde n]$ of
Laurent polynomials in $q$. Hence the action of $U'_q\so_n$
on its Verma module $V_m$ is also well-defined for $q$ a root of unity,
with the coordinates of $m$ and $\tilde n$ being complex numbers.
 $\Box$

\medskip
To give some indication about the additional phenomena at roots
of unity, we prove a special case
of the main result of \cite{IK2} in our setting.

\begin{theorem}\label{lgeneric} Let $\la=(\la_i)$ be a weight such that
$\la_i\not\in \frac{1}{4}\Z$, $1\leq i\leq n/2$,
and let $m=[\la]$ or $m=[\la]_+$. Then there exists a
module $M_m$ of dimension $\ell^d$,
where $d$ is the number of positive roots
in $\so_n$ and with highest weight $\la$.
\end{theorem}

$Proof.$ The proof of this theorem is a variation of the one of Theorem
\ref{character}. First of all observe that all the coefficients are
well-defined in the explicit representations of $\so_3$ and $\so_4$
in this and the previous sections. One directly reads off these representations
that $v_{\la-\ell\al_i}$ is a well-defined highest weight vector in
the Verma module $V_m$ for $\so_4$ and $\so_3$, a result analogous to
the one of Proposition \ref{Vermaso4}. For $n>4$, we can now define
a submodule $V_i$ with highest weight $\la-\ell\al_i$ as in Lemma
\ref{generalideals}.
We now define $M_m=V_m/\sum V_i$.

To determine the dimension and character of $M_m$, we map the weights
of $V_\la$ to the weights of $V_{(\ell-1)\rho}$, via the map
$\mu\mapsto \mu-\la+(\ell-1)\rho$. Using this map, we
pull back the dot-action of the Weyl
group on the weight lattice
to the lattice $\mu +P$, where $P$ is the weight lattice.
Then one shows as in the proof of
Theorem \ref{character} that the weight spaces of $M_m$ are symmetric
under this action of the Weyl group. In particular,
we obtain that
$$\dim M_m = \dim V_{(\ell-1)\rho} = \prod_{\al>0}\frac{(\ell\rho,\al)}{(\rho,\al)}\ = \ell^d.$$

\begin{remark} The same proof also works if we define $V_i(a_i)$ to be
the Verma module generated by $a_iv_\la + v_{\la-\ell \al_i}$, where
$v_{\la-\ell \al_i}$ is the highest weight vector of the module $V_i$
in the proof of Theorem \ref{lgeneric} and $a_i\in\C$. Hence we obtain a
multiparameter family of modules $V/\sum_i V_i(a_i)$
with highest  weight $\la$. It is easy to see
that these modules do not have a lowest weight vector if $a_i\neq 0$
for some $i$. So we have many non-isomorphic finite-dimensional 
modules with the same highest weight and the same dimensions.

There are even more non-isomorphic modules with the same weight structure,
see \cite{IK2}, which do not have highest weight vectors.
These seem to be phenomena quite parallel to the representation theory
of quantum groups at roots of unity. E.g.  the modules in Theorem \ref{lgeneric}
appear to be analogs of what are called baby Verma modules 
for quantum groups. 
\end{remark}

\subsection{Unitary representations} One of the motivations for
this paper came from the need to identify certain representations
of $U'_q\so_n$ on a Hilbert space.
Fortunately, in spite of all the new phenomena we have seen in the last subsection,
the situation for unitary representations will turn out to be similarly nice as for $q$ not a root
of unity.  See Example 2 in Section
\ref{knownrep} for the choice of name in the following definition.
We call a representation of $U'_q\so_n$
on a Hilbert space $V$ with inner product $(\ ,\ )$ a {\it unitary}
representation if
\begin{equation}\label{unitform}
(B_iv,w)=(v,B_iw)\quad {\rm for}\  v,w,\in V,\ 1\leq i\leq n-1\quad {\rm for \ all\ } v,w\in V.
\end{equation}
We are going to show that unitary representations with highest
weights will have to factor
over an ideal like $I(\la)$ in Theorem \ref{character}. The idea is to prove
that the image of the highest weight vectors $v_{s_i.\la}$ in a unitary
representation will have to have length 0, using the explicit representations for $U'_q\so_3$ and $U'_q\so_4$. There are some minor complications
as the weight vectors there are not always well-defined at a root of unity.

\begin{lemma}\label{exclude} Let $(\ ,\ )$ be a bilinear or sesquilinear form
on the Verma module $V_m$ of $U'_q\so_3$ with respect to which the generators 
$B_1$ and $B_2$ are self-adjoint.

(a) If the weight vectors $v_i$ are well-defined and mutually orthogonal for $j-1\leq i\leq j+2$, then $\| v_{j+1}\|^2=\al_{j,j+1}\| v_j\|^2$.

(b) The weight $m=[\pm \ell/4]=\pm[0]_+$
can appear in a unitary highest weight representation $V$ of $U'_q\so_3$
only for highest/lowest weight vectors. In particular, if $w_j\in V$ is a weight vector
with weight $\pm [0]_+$, then $B_2w_j$ is a weight vector with weight $\pm [1]_+$
and with $\| B_2w_j\|^2=-2/(q-q^{-1})^2\| w_j\|^2$.

(c) If $V$ is a unitary representation of $U'_q\so_3$ with highest weight $m=[\la]$,
and $\Phi: V_m\to V$ the canonical epimorphism, then $\Phi(v_{2\la +1})=0$.
\end{lemma}

$Proof.$  Using the assumptions and the definition of the action of $B_2$, 
see Proposition \ref{weightbasis}, we obtain
\begin{equation}\label{lengths}
(v_{j+1},v_{j+1})=(B_2v_j-\al_{j-1,j}v_{j-1},v_{j+1})
=(v_j, B_2v_{j+1})=\al_{j,j+1}(v_j,v_j).
\end{equation}
To prove (b),
assume to the contrary that a weight $\pm [0]_+$ does
occur in a module with highest weight $[\la]_+$
for the weight vector $v_j$, $0<j$ such that $B_2v_j$
is not a multiple of $v_{j-1}$. Observe that $[\la-j]_+=[0]_+$
implies $j=\la$.
We proceed as in Lemma \ref{welldefined} by
defining the vector $v'_{j+1}=B_2v_j$. One calculates from the
definitions that
\begin{align}
B_1v'_{j+1}\ &=\ [\la-j-1]_+v'_{j+1}+(q-q^{-1})[\la-j]\al_{j-1,j}v_{j-1}\ =\cr
&=\ [\la-j-1]_+v'_{j+1}+ i\frac{[\la+1][\la]}{q-q^{-1}}v_{j-1},
\end{align}
where we only used $\la=j$ for the second equality.
Observe that $[\la-j-1]_+=[\la-j+1]_+$. If $[\la+1][\la]\neq 0$,
$B_1$ would act as
a Jordan block on span $\{ v'_{j+1}, v_{j- 1}\}$, a contradiction to it being
self-adjoint. If $[\la+1][\la]=0$, we can show that a vector of weight $[0]_+$
has to be equal to 0 as follows:

If $[\la+1]=0$,  consider $m=[r\ell/2-1]_+$, where $\ell$ is even and $r$ is a positive integer.
One checks, using the formulas of $\al_{j-1,j}$ in Corollary \ref{ali}
that the vectors $v_j$ are well-defined for
$\ell/2-3\leq j\leq \ell/2$. One deduces from this that $\|v_{\ell/2-1}\|=0$,
as $\al_{\ell/2-2,\ell/2-1}=0$. Hence a unitary module with highest weight
$[(2r-1)\ell/4-1]=[r\ell/2-1]_+$ only contains the weights $[r\ell/2-j]_+$ with $1\leq j\leq \ell/2-1$.
 If $\ell$ is odd and $[\la+1]=0$,
$\la$ would have to be of the form $\la=r\ell-1$, for some integer $r>0$. But then already
$\al_{0,1}=0$, which implies $\|v_1\|=0$.

For the case $[\la]=0$: If $m=[\ell/2]_+=[\ell/4]$ with $\ell$ even, one checks, using the formula for
$\al_{\ell/2-1,\ell/2}$ with $m=[\ell/4]$ that
$v_{\ell/2+1}= B_2v_{\ell/2}-\frac{-2}{(q-q^{-1})^2}v_{\ell/2-1}$
is a well-defined weight vector. Using the recursion formula
\ref{recursealpha} for $j=\ell/2$ and $m_{j}=\pm[0]_+$, we see
that $\al_{\ell/2,\ell/2+1}=0$ also for our choice of $q$. Unfortunately, we can not conclude
that $\|v_{\ell/2+1}\|=0$, as $v_{\ell/2-1}$ need not necessarily be orthogonal
to $v_{\ell/2+1}$, having the same weight. But we can conclude that it would also 
be a highest weight vector. By the previous results, the module $U$ generated by such
a vector in a unitary highest weight representation $V$ would not contain the weight $\pm[0]_+$.
Hence the image of the highest weight vector $v_0$ would be in the orthogonal complement of $U$.
This would contradict the fact that $v_0$ generates $V$, unless $U=0$.
 Hence the image of the vector $v_{\ell/2}$
in a unitary representation would be a highest/lowest weight vector, and it would be the only weight vector
besides $v_0$ with weight $\pm [0]_+$.
The norm of the weight vector $B_2w_{\ell/2}$ can now be calculated
as in (a).

For part (c) we can calculate the norms of the weight vectors $v_j$ by
(a), as long as the weights $\pm [0]_+$ do not appear. By (b), this
is the case for all unitary highest weight representations except
when the highest weight is equal to $\pm [0]_+=[\pm \ell/4]$.
If $m=[\la]\neq [\pm \ell/4]$, it follows as in the classical case that
$\|v_{2\la+1}\|^2=\al_{2\la,2\la+1}\|v_{2\la}\|^2=0$. If $m=[\pm \ell/4]$,
then $\Phi(v_{\ell/2})$ is a lowest weight vector  and 
$\Phi(v_{\ell/2+1})=\Phi(v_{\ell/2})-\frac{2}{(q-q^{-1})^2}\Phi(v_{\ell/2-1})=0$,
by (b). $\Box$

\begin{lemma}\label{so4unity} Let $m=[\la]$ 
or $m=[\la]_+$ with $\la$ as in Theorem \ref{compclass}
additionally 
satisfying $\ell/4\geq \la_1\geq\ ...\ \geq
|\la_{\lfloor n/2\rfloor}|$.
Then $\| v_{s_i.\la}\|=0$, where $v_{s_i.\la}$
is as in Proposition \ref{Vermaso4}.
\end{lemma}

$Proof.$
Let $m=[\la]$,
let $\mu=\la-r\al_1$, and let $v_\mu$ be the weight vector as
defined in Lemma \ref{so4normalization}, (b). Observe that these vectors are 
well-defined also for $q$ a root of unity.
We are going to use Eq \ref{so4simple} for the action of $B_2$ on $v_\mu$.
Observe that even if $\{\mu_2\}=0$, the expression
$$\frac{1}{\{\mu_2\}}(v_{\mu-\al_1}-v_{\mu-\al_2})= B_2v_\mu -\frac{[r_1][(\la,\al_1)+1-r]}{\{\mu_1\}}v_{\mu+\al_1}$$
is well-defined, provided $\{\mu_1\}=[\mu_1+1]-[\mu_1-1]\neq 0$. 
In particular, in this
case this expression is an eigenvector of $B_1$ with eigenvalue $[\mu_1-1]$,
and hence orthogonal to $v_{\mu+\al_1}$.
Using this for the last equality below, as well as Eq. \ref{so4simple}, we obtain
\begin{align}
\| v_{\mu-\al_1}\|^2\ &=\ ((B_3-[\la_2+r-1])B_2v_\mu,v_{\mu-\al_1})=
(v_\mu, B_2(B_3-[\la_2+r-1])v_{\mu-\al_1})\ =\cr
&=\{\la_2+r\}(v_{\mu}, B_2v_{\mu-\al_1})\ =\cr 
&=\frac{\{\la_2+r\}}{\{\la_1-r-1\}}[r+1][\la_1-\la_2-r]
\ \| v_\mu\|^2.
\end{align}
Now observe that if 
$r=(\la,\al_1)=\la_1-\la_2$ and $\{\la_1-r-1\}=\{\la_2-1\}\neq 0$,
we have $s_1.\la=\mu-\al_1$ and hence $\|v_{s_1.\la}\|^2=0$.
So except if $\la_2=1-\ell/4$, we have shown that
$v_{s_1.\la}$ is in the kernel of $\Phi$.
One shows by the same approach
that also $\|v_{s_2.\la}\|^2=0$ as long as $\la_2\neq \ell/4-1$.
This applies, in particular, to the remaining cases with
$\la_2=1-\ell/4$. We leave it to the reader to check that
the quotient $V_\la/V_{s_2.\la}$ is isomorphic to $V_{\ell/4-1}$
as a $U'_q\so_3$-module for $\la_1=\ell/4-1$, and it is isomorphic to
$V_{\ell/4}+V_{\ell/4-1}$ as a $U'_q\so_3$-module for $\la_1=\ell/4$.
One can now check for $\la=(\ell/4,\ell/4-1)$ that the vector
$v_0=[(\la,\al_2)]v_{\la-\al_1}+[(\la,\al_1)]v_{\la-\al_2}$
does generate the $U'_q\so_3$ Verma module with highest weight $[\la_1-1]$,
and, using Lemma \ref{exclude}, that the vector $v_{\ell/2-1}$
in this Verma module is a nonzero multiple
of $v_{s_1.\la}$ mod $V_{s_2.\la}$. Lemma \ref{exclude} then implies
that $\| v_{s_1.\la}\|=0$. The case with $\la=(\ell/4-1,1-\ell/4)$ 
is similar and easier.  $\Box$

\begin{theorem}\label{factorweyl} Let $m=[\la]$ 
or $m=[\la]_+$ with $\la$ as in Theorem \ref{compclass}
additionally 
satisfying $\ell/4\geq \la_1\geq\ ...\ \geq
|\la_{\lfloor n/2\rfloor}|$.
Then there exists at most one simple unitary $U'_q\so_n$ module
for $m=[\la]$ and at most $2^{\lfloor (n-1)/2\rfloor}$ 
simple unitary modules  for $m=[\la]_+$
which would be a quotient of the standard 
Verma module $V_m$.
\end{theorem}

$Proof.$ Let $\Phi$ be the surjective map from the Verma module $V_m$ onto
a unitary simple representation $V$ with
highest weight $m$, and let $(\ ,\ )$ also denote the pull-back
of the inner product on $V$ to $V_m$. Using subalgebras
$U'_q\so_4(i)$ and $U'_q\so_3(i)$ as in the proof of Lemma
\ref{generalideals} and Theorem \ref{character},
we conclude from Lemmas \ref{so4unity} and \ref{exclude}
that $\| v_{s_i.m}\|^2=0$ for the vectors $v_{s_i.m}$ defined there.
Hence the ideal $I(\lambda)$ defined in Theorem \ref{character}
must be in the kernel of $\Phi$. The character of $V_m/I(\lambda)$ is 
given by Weyl's character formula. If $m=[\la]$, any weight
$[\mu]$ satisfies $|\mu_i|\leq \ell/4$, by our assumption on $\lambda$.
So even for $q$ a primitive $\ell$-th root of unity, we have
$[\mu]=[\nu]$ for two weights $\mu$ and $\nu$ of $V$
if and only if $\mu=\nu$. It follows that the weight space of the highest weight $[\la]$
is one-dimensional.  We can now use the standard 
argument for showing that there exists a unique maximum submodule $J$ of $V_m$
with $I(\la)\subset  J$. Indeed, if there were two such submodules $J_1$ and $J_2$,
then also $(J_1+J_2)/I(\la)$ would not contain the highest weight, and hence
also $J_1+J_2$ would be a maximum submodule. As $V$ is simple,
the kernel of $\Phi$ has to be equal to the maximum submodule $J$.
This implies the claim for $m=[\la]$. If $m=[\la]_+$, $V_m/I(\la)$
splits into the direct sum of $2^{\lfloor (n-1)/2\rfloor}$ modules,
each with highest weight $m$ (see Corollary \ref{Weylmplus}).
We can similarly
show that the highest weight $[\la]_+$ appears only with multiplicity one
in any of these summands. One now shows as before that any maximum
submodule $J\supset I(\la)$ has to contain all but one of these
summands and the unique maximum submodule of the remaining summand.
 $\Box$

\begin{remark} 1. It can be shown that there are indeed unitary representations
of $U'_q\so_n$ for the highest weights $m$ listed in Theorem \ref{factorweyl} for
 suitable $\ell$-th roots of unity of $q$. If $m=[\la]$, these would be the $\ell$-th roots  closest 
to 1, while for $m=[\la]_+$ they would be the ones closest to $-1$. In the first case,
this follows from the representations mentioned under 2. in Section \ref{knownrep},
using the results in \cite{Wcat}. As mentioned there, analogous non-classical
representations have only been found recently and are currently being worked out in \cite{Wdual}.
Assuming this result, one can also deduce the existence of unitary non-classical representations
as before, using \cite{Wcat}. It may also be possible to give a more direct proof
closer to the spirit in this paper by proving that the bilinear
form defined in \ref{unitform} is positive semi-definite on $V_m$ for suitable values of $q$.

2. This last theorem can be used to reprove the main technical
result of \cite{RW}. It was shown there that the unitary representations
of $U'_q\so_n$ of Section \ref{knownrep}, Example 3 had the same characters, up to rescaling multiplicities, as certain representations listed under Example 2.
Hence they are isomorphic.

3. For a root of unity, the representations of Theorem \ref{character} are usually not simple. Assuming that the duality between the actions of $U_q\so_N$ and $U'_q\so_n$ on $S^{\otimes n}$ also holds for $q$ a root of unity,
where $S$ is the spin representation of $U_q\so_N$ (see Example 2 and \cite{WSp} for details), the multiplicities of
simple modules in a filtration should be given by certain parabolic Kazhdan-Lusztig polynomials. The latter ones have appeared to describe characters of tilting modules of Drinfeld-Jimbo quantum groups.
\end{remark}

\section{Discussion of our results}\label{discussresults}

\subsection{Technicalities in this approach} The main results
in this paper concern the construction and classification of
highest weight modules of $U'_q\so_n$ via a
Verma module approach. One of the difficulties in this approach
is the lack of
generic raising and lowering operators: One can define, for
a given weight vector $v$ of weight $\mu$ and a given simple root $\alpha$,
an expression $E$ in terms of our generators such that $Ev$ has weight
$\mu+\alpha$ (if it exists), see \cite{Kl}. 
But this expression depends on the given weight $\mu$;
the same expression applied to another weight vector usually does
not even give a weight vector anymore. This has made our approach more
involved than originally expected. In particular, we have decided to appeal
to already known representations e.g. for proving linear independence
of our spanning set in general. It would be interesting to see whether
this could be avoided without too much additional work. On the other hand,
finding sufficiently many representations of 
general coideal subalgebras of quantum groups should
not be terribly hard in view of the established representation theory
of quantum groups; we only needed the classical representations
for proving linear independence.
Another short-coming in our approach is the fact
that we did not prove irreducibility of the finite-dimensional modules
coming from our quotient construction. This should not be overly hard to fix
using Casimir elements for $U'_q\so_n$, see Remark \ref{irreducibilityrem} for details.
For the cases treated here, it would also follow from the results in
\cite{WSp} and \cite{Wdual}, besides the original result in \cite{IK}.

\subsection{Comparison with the approach by Klimyk et al} Recall our
brief discussion in Section \ref{defrep}
of the classification of representations
of $U'_q\so_n$ in \cite{IK}. We think the approach in this paper
has two advantages. It also works for representations at roots
of unity with integral dominant weight, allowing us to identify
at least certain representations at roots of unity via their highest
weight. These include the ones studied in \cite{RW}.
Secondly, the approach in this paper is less
computational. While this may just be a matter of personal taste,
we hope that the approach here may be more suitable
for generalizations to other coideal algebras (see also next
subsection). On the other hand, one can criticize our approach
for not proving irreducibility of the finite-dimensional modules.
But see the remark at the end of the previous subsection and
Remark \ref{irreducibilityrem}.

\subsection{Generalizations for other coideal algebras}
More general coideal subalgebras of quantum groups appear
in several contexts where it would be useful to know their
representation theory such as e.g. in $q$-Howe duality
and categorification of representations of quantum groups
(see e.g. \cite{ES}, \cite{ST}). So it is a natural
problem  to classify the representations of general coideal
subalgebras.
As a first general result, Letzter has determined analogs of Cartan
subalgebras for all coideal subalgebras in \cite{LC}.
It would be interesting to see whether the approach
in this paper could be extended to a more general setting.

\end{document}